\documentclass[11pt]{article}
\usepackage[left=1in,right=1in,top=1in,bottom=1in]{geometry}

\usepackage[normalem]{ulem} 
\usepackage[utf8]{inputenc} 
\usepackage[T1]{fontenc}    
\usepackage{hyperref}       
\usepackage{url}            
\usepackage{booktabs}       
\usepackage{nicefrac}       
\usepackage{microtype}      

\usepackage{graphicx, xcolor, ulem}
\usepackage{multirow}
\usepackage{enumitem}
\usepackage{amsmath,amsfonts,latexsym,amssymb,amsbsy, amsthm,textgreek,dsfont}

\usepackage{footnote}
\usepackage{thmtools, thm-restate}
\usepackage{authblk}

\usepackage{url, color, xcolor, amsmath}
\usepackage[customcolors,shade]{hf-tikz}

\DeclareMathOperator*{\argmax}{\arg\!\max}
\newcommand{\E}{\mathop{\mathbb{E}}}
\newcommand{\PR}{\mathop{\mathbb{P}}}
\renewcommand{\deg}{\mathop{\mathrm{deg}}}
\usepackage{algorithm}
\usepackage[noend]{algpseudocode}
\algdef{SE}[SUBALG]{Indent}{EndIndent}{}{\algorithmicend\ }%
\algtext*{Indent}
\algtext*{EndIndent}

\newcommand*{\QED}{\null\nobreak\hfill\ensuremath{\square}}

\newtheorem{theorem}{Theorem}

\newtheorem{lemma}{Lemma}

\newtheorem{definition}{Definition}

\algnewcommand\INPUT{\item[{\textbf{Input:}}]}
\algnewcommand\RETURN{\item[{\textbf{Return:}}]}
\newcommand\numberthis{\addtocounter{equation}{1}\tag{\theequation}}

\allowdisplaybreaks 
\title{Fair and Reliable Reconnections for Temporary Disruptions in\\ Electric Distribution Networks using Submodularity}
\author[1]{Cyrus Hettle}
\author[1]{Swati Gupta}
\author[1]{Daniel Molzahn}

\affil[1]{\small Georgia Institute of Technology \protect \\
{\small \tt \{chettle,swatig,molzahn\}@gatech.edu}}
\date{}

\begin{document}
\maketitle
\begin{abstract}Increasing reliability and reducing disruptions in supply networks are of increasing importance; for example, power outages in electricity distribution networks cost \$35-50 billion annually in the US. Motivated by the operational constraints of such networks and their rapid adoption of decentralized paradigms and self-healing components, we introduce the “minimum reconnection time” (MRT) problem. MRT seeks to reduce outage time after network disruptions by programming reconnection times of different edges (i.e., switches), while ensuring that the operating network is acyclic. 

We show that MRT is NP-hard and is a special case of the well-known minimum linear ordering problem (MLOP) in the submodular optimization literature. MLOP is a special case of a broader class of ordering problems that often admit polynomial time approximation algorithms. We develop the theory of kernel-based randomized rounding approaches to give a tight polynomial-time approximation for MRT, improving the state-of-the-art approximation factor for a broad class of MLOP instances. Further, motivated by the reliability incentive structure for utility companies and operational energy losses in distribution networks, we propose local search over spanning trees to balance multiple objectives simultaneously. We computationally validate our reconfiguration methods on the NREL SMART-DS Greensboro synthetic network, and show that this improves service equity by a factor of four, across industrial and residential areas.
\end{abstract}

\section{Introduction}

Improvements in distribution network reconfiguration are motivated by the high societal cost of electricity interruptions. Despite investments from utility companies in improving reliability, service disruptions remain challenging~\cite{larsen2020}, 
with the economic costs of power outages estimated between \$35 billion and \$50 billion annually in the United States~\cite{lacommare2018}. Related studies on the Value of Lost Load (VOLL) give per-kilowatt-hour estimates ranging from several dollars to several hundred dollars, depending on the underlying methodology, customer type, and assumptions, and the Midcontinent Independent System Operator, one of the largest power system operators in the United States, uses a VOLL of \$3.50 per kilowatt-hour~\cite{schroder2015,miso_voll}. These values are an order of magnitude larger than typical costs (around 10.5 cents per kilowatt-hour in the United States~\cite{eai2019}), reflecting the importance of providing reliable electricity.

Since reliability is so valuable, industry groups recognize utilities using reliability metrics such as the System Average Interruption Duration Index (SAIDI), e.g., the ReliabilityOne Award and the American Public Power Association's Certificate of Excellence Award~\cite{reliabilityone,appa_awards}. Some regulators implement explicit financial incentives. For example, the California Public Utilities Commission set a SAIDI target of 60 minutes for San Diego Gas \& Electric~\cite{sdge_rates}. Overperforming or underperforming this target beyond $\pm 2$ minutes was rewarded or penalized by \$375,000 per minute up to a maximum of \$$3$ million. Similar incentives exist in Great Britain, Italy, etc. ~\cite{regulation_survey,NBERc12566}. 

Reducing outage time is an important way to improve a network's reliability metrics. After an isolated fault occurs in a system, causing an outage downstream in the network, it is crucial to quickly reconfigure the network to reroute power to the affected areas. Reconfiguration can potentially be executed long before the fault is manually repaired, reducing the outage duration for many customers. Much of the existing work in the academic literature focuses on \emph{centralized} approaches for distribution network reconfiguration. In centralized approaches, information regarding distribution system operations, such as the location of a fault that is causing a power outage, is communicated to a central substation or control center which computes an appropriate switching configuration that is then communicated to each switch. In contrast, \emph{decentralized} approaches have local controllers that operate solely based on local information. Solving a reconfiguration problem centrally with full information and control is simpler, but has speed, robustness, and security disadvantages \cite{wang2017distributed,kroposki2020autonomous}. For these reasons, protection and reconfiguration systems are moving towards an adaptive paradigm where components of the distribution system work together to detect and characterize failures and then restore the supply of power as quickly as possible, e.g., ~\cite{mishra2017,sakis2017,vanderwalt2018}.%

To this end, we propose a decentralized approach for automatically reconfiguring radial distribution networks to quickly restore power after a fault. Such models are classified as Fault Location, Isolation, and Service Restoration applications that embed ``self-healing'' features within Advanced Distribution Management Systems~\cite{epri_onr}. 
A variety of self-healing approaches use decentralized schemes; e.g., \cite{li2020} isolate faults onto ``islands'' and use alternate paths to restore service, and ~\cite{torres2018distributed} employ an internal timer for each switch, as in our model. We assume a network is given as a graph $G = (V,E)$ with an initial radial tree $T \subseteq E$. We employ edges (e.g., switches in an electrical network) $s\in S = E\setminus T$ that are not in the base radial network, each of which can detect abnormal flow patterns on its endpoints, and then activate to reconnect the network. Our decentralized approach operates by preassigning an ordering to the edges in $S$ that describes how long each edge waits before attempting to activate. This formulation leads to a novel connection of network reliability to linear ordering problems from the submodular optimization literature. 

\begin{figure}[!t]
 \begin{center}
    \newcommand{\scale}{.5} 
\begin{tikzpicture}[scale=\scale]
         \draw[black](1,0)--(2.5,2);
         \draw[black](4,0)--(2.5,2);
         \draw[black](0,2)--(2.5,2);
         \draw[black](5,2)--(2.5,2);
         \draw[black](1,4)--(2.5,2);
         \draw[black](4,4)--(2.5,2);
        \draw[black,dashed](1,0)--(4,0);
         \draw[black,dashed](4,0)--(5,2);
         \draw[black,dashed](0,2)--(1,0);
         \draw[black,dashed](5,2)--(4,4);
         \draw[black,dashed](1,4)--(0,2);
         \draw[black,dashed](4,4)--(1,4);
        \node at (1,4) {$\bullet$};
         \node at (4,4) {$\bullet$};
       \node at (1,0) {$\bullet$};
       \node at (4,0) {$\bullet$};
       \node at (0,2) {$\bullet$};
       \node at (2.5,2) {\huge$\bullet$};
       \node at (5,2) {$\bullet$};
       
        \draw[black,dashed ](9,0)--(10.5,2);
         \draw[black,dashed ](12,0)--(10.5,2);
         \draw[black,dashed ](8,2)--(10.5,2);
         \draw[black,dashed ](13,2)--(10.5,2);
         \draw[black,dashed ](9,4)--(10.5,2);
         \draw[black ](12,4)--(10.5,2);
        \draw[black ](9,0)--(12,0);
         \draw[black ](12,0)--(13,2);
         \draw[black ](8,2)--(9,0);
         \draw[black ](13,2)--(12,4);
         \draw[black ](9,4)--(8,2);
         \draw[black,dashed ](12,4)--(9,4);
        \node at (9,4) {$\bullet$};
         \node at (12,4) {$\bullet$};
       \node at (9,0) {$\bullet$};
       \node at (12,0) {$\bullet$};
       \node at (8,2) {$\bullet$};
       \node at (10.5,2) {\huge$\bullet$};
       \node at (13,2) {$\bullet$};

    \end{tikzpicture}
            
\end{center}
\caption{Two spanning trees on the same network. Solid edges are in $T$, dashed edges are in $E\setminus T$, the largest vertex is  source $r$, and all vertex weights and edge failure probabilities are equal.}
\label{fig:spokewheel}
\end{figure}
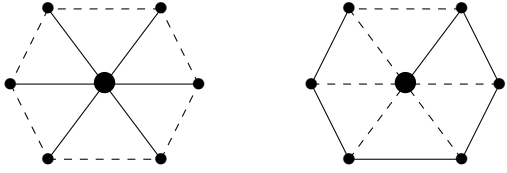

\paragraph{Contributions.} 
We summarize our contributions next. In section~\ref{sec:model}, we present our general graphical model for network flow and network disruptions. We give a novel system for reconnecting the network after a disruptive event by activating edges in a decentralized manner, according to a preprogrammed order. This model is based on the challenges seen in fault detection and protection systems in electricity distribution networks. We define the Minimum Reconnection Time (MRT) problem of finding the best activation order of edges. Different fine-based structures imposed on utility companies, such as SAIDI, can be written as objectives of instances of MRT.

We show in section~\ref{sec:mssc} that MRT is in fact NP-hard, using a reduction from the {\it tree augmentation problem} (TAP) of selecting a minimum-cardinality subset $F$ of a set of edges $E$ such that $T\cup F$ is 2-edge-connected, where $T$ is a given tree~\cite{Cohen2013}. Next, we show that the MRT problem is connected to the field of submodular optimization (e.g., \cite{lovasz1983submodular,edmonds1970submodular}), and in particular is a special case of the Min Sum Set Cover (MSSC) problem~\cite{Feige2002}. Submodular functions are a discrete analog of convexity and capture the property of diminishing returns, and so are widely used in applications such as scheduling~\cite{happach2020min} and network topology~\cite{narayanan1997submodular}. 

In section~\ref{sec:approx}, we give an improved approximation factor of $(2c/(c+1))^2$ for MRT instances which are {\it $c$-uniform}, i.e., for each initial radial tree edge $e\in T$, there are at most $c$ non-tree edges $s\in S=E\setminus T$ such that $T-e+s$ is connected. This improves on the best known approximation factor of 4 for $c$-uniform instances \cite{Feige2002}. To achieve this improvement, we use {\it fractional kernels} within a kernel-based randomized rounding approach~\cite{phillips1998minimizing,hall1997scheduling}, which may be of independent interest. These kernels are based on discrete analogs of negative power functions, and are used to spread the weight of the relaxed optimal solutions across time. Since our kernel is fractional, the analysis of our method is more intricate and uses discrete calculus-like manipulations. We further show that the $(2c/(c+1))^2$ approximation ratio is tight (Theorem~\ref{thm:alphatight}) for all $c$-uniform MSSC instances, with respect to the natural LP relaxation.   

Given a network, the choice of tree $T$ can have a large impact on the reconnection metrics and losses due to energy dissipation. For example, Figure~\ref{fig:spokewheel} shows two different spanning trees on the 7-vertex wheel graph. The ``spoke'' network (left) achieves an objective value $\frac{3}{2}$ for both {\sc R-Time} (an outage metric we introduce) and SAIDI, while the ``wheel'' network (right) achieves objective value $1$ for {\sc R-Time} and $\frac{7}{2}$ for SAIDI. Moreover, the choice of the initial tree determines the energy loss in the network, which can be as large as 5-10\% of the total energy generation~\cite{distribution_energy_losses,distribution_energy_losses_world_bank}.
Energy loss is 6 in the ``spoke'' and 91 in the ``wheel,'' where we assume each bus has a demand of one unit, lines have uniform resistance, and losses are proportional to the square of the power flow through each line. Moreover, we find that on the NREL SMART-DS Greensboro synthetic network, optimizing for {\sc R-Time} induces an expected outage time for medium voltage buses that is 50\% longer than when optimizing for SAIDI. Therefore, even the choice of objectives can have a non-negligible disparate impact on service for different ``groups'' in the network.

Therefore, to optimize multiple metrics simultaneously, in section~\ref{sec:local} we use local search techniques over the set of spanning trees with a multi-criteria objective that involves all three metrics of interest (SAIDI, {\sc R-Time}, and energy losses). We show that the choice of optimization metric results in varying outage times for consumers at low and medium-voltage buses. Finally, in section \ref{sec:computations}, we test our methods on the NREL SMART-DS Greensboro synthetic network. By adding relatively few edges to the network (less than .2\% of the number in the original network), prioritizing the heaviest-load areas, we show that we can 
improve connectivity to cover over 90\% of load in the network. After adding edges, our local search methods significantly improve outage and energy metrics, with up to a 95\% reduction in the product of {\sc R-Time}, SAIDI, and energy loss and a nearly 60\% reduction on average. As a byproduct, we show that the equity of service between industrial and residential areas (i.e., medium and low voltage buses) improves by a factor of at least 4. We believe that next-generation networks with self-healing capabilities will benefit from decentralized approaches, as well as multi-criteria optimization to balance various reliability metrics.

\section{Problem formulation and modeling}\label{sec:model}
We state our model in terms of general networks, where flow on an edge might represent electricity, gas, water, goods, etc.
We give additional motivation and modeling details for electrical networks in Appendix~\ref{app:elec_networks}.

Consider a network represented as a graph $G=(V,E)$, with vertex set $V$ and edges $e\in E$. Each edge can be {\it active} (allows flow) or {\it inactive} (no flow allowed, e.g., an open switch). We assume that there is a given initial configuration of active edges $T \subseteq E$ which is a spanning tree on $V$ and meets the flow demands in the network. Let $r \in V$ be a designated root vertex (e.g., an electrical substation) which supplies the demand in the network to the other vertices. Initially, the network operates using edges in $T$, but once a disconnecting ``event'' occurs, e.g., a fault in an electricity distribution network caused by an isolated event \cite{larsen2015}, an edge goes from being active to being inactive. As a result, the set of active edges becomes a disconnected union of two trees.

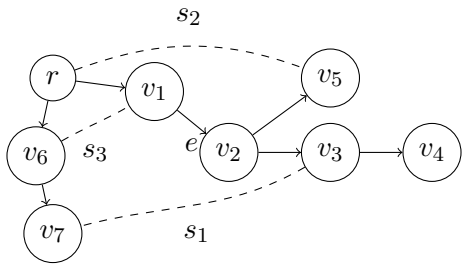
\begin{figure}[t]
\begin{center}
\begin{tikzpicture}[scale=0.45]
\tikzset{every node/.style={draw=black,circle}};
    \node (s) at (-3,0) {$r$};
    \node (6) at (0,-.4) {$v_1$};
    \node (7) at (2.2,-2.2) {$v_2$};
    \node (8) at (5.2,-2.2) {$v_3$};
    \node (9) at (8.2,-2.2) {$v_4$};
    \node (10) at (5.2,0) {$v_{5}$};
    \node (11) at (-3.5,-2.3) {$v_{6}$};
    \node (12) at (-3,-4.6) {$v_{7}$};
    \draw[->] (7) to (8);
    \draw[->] (8) to (9);
        \draw[->] (s) to (11);
        \draw[->] (s) to (6);
    \draw[->] (11) to (12);
    \draw[->] (7) to (10);
    \tikzset{switch/.style={black,fill=none,draw=none,midway}};
    \draw[->] (6) to  node[switch,below] {$e$} (7);
    \draw[dashed] (8) to[out = -150, in = 15] node[switch,below] {$s_1$} (12);
    \draw[dashed] (s) to[out = 20, in = 160] node[switch,above]{$s_2$} (10);
      \draw[dashed] (11) to[out = 30, in = 210] node[switch,below]{$s_3$} (6);
\end{tikzpicture}
    \end{center}
    \caption{An example of a network comprising a spanning tree (solid edges) and inactive edges (dashed edges).}
    \label{fig:example_network_1}
\end{figure}
 
When an event occurs at edge $e \in T$, causing it to become inactive, all vertices that receive flow via $e$ will become disconnected from the root vertex $r$. Let $R\subseteq V$ be the {\it reachable set} of vertices, i.e., those in the same connected component as $r$ in the subgraph of $G$ comprising the active edges.  Activating an edge $s \in S=E\setminus T$ could allow flow through $s$ from $R$ to $V \setminus R$. This is feasible if and only if $e$ is in the unique cycle of $T + s$; in other words, $T - e + s$ is a tree. In this case, we say that $s$ \textit{covers} $e$. To maximize service availability (for instance, before potentially time-consuming repairs necessary to reactivate $e$ are completed), we would like to activate such an edge to restore flow as quickly as possible. We must do so while respecting the requirement that the connected component containing $r$ is radial (i.e., has no cycles), which is an operational requirement often seen in electricity and gas networks~\cite{dokic_integrated}.

\textit{Example.} Consider the network in Figure~\ref{fig:example_network_1}. Vertex $r$ is the root. Edges in $T$ are solid and edges in $S$ are dashed. Suppose an event occurs on edge $e$, deactivating it. Vertices $v_2,v_3,v_4,$ and $v_5$ will be without flow since they are no longer in the same component as the root $r$. If $s_1$ or $s_2$ are activated, then the flow would be restored. Activating $s_3$ would not restore flow, and activating both $s_1$ and $s_2$ simultaneously would violate the radial constraint. 

In our model, choosing an edge to activate would be straightforward given perfect information on the status of the network and the ability to activate edges through a central planner. However, there is an increasing push in many settings for decentralized operations, where cheap components such as sensors can detect local changes to the operability of the network (e.g., an edge can detect whether it has flow going through itself). Based on this information, an automatic self-healing operation can be performed, such as activating an edge. Therefore, we consider the following localized and decentralized reconnection process. Each non-tree edge $s\in E\setminus T$ has a sensor which can detect whether each of its incident vertices is receiving flow (i.e., is in the reachable set $R$). Upon detecting that exactly one of its incident vertices is in $R$, it waits for a preprogrammed number $\sigma(s)$ of time steps. If during that period of time that status of its incident vertices with respect to $R$ does not change, the edge activates itself.

To define the preprogrammed waiting times for each {\it self-healing} edge, we select a permutation $\sigma:S\longrightarrow [|S|]$, where $[|S|]=\{1,2,\dots,|S|\}$. If necessary, we choose time steps to be long enough that any changes in the network due to other edges activating, flow being restored, etc. will be detected by the sensors within one unit of time. This ensures that multiple edges do not activate in such quick succession that a cycle forms in the set of active edges. Note that any ordering $\sigma$ will yield feasible values for $\sigma(s)$ and eventually reconnect the network and restore flow safely if it is possible to do so (i.e., without creating a cycle), but we would like to do so as quickly as possible so as to minimize network reconnection fines and other operational objectives such as energy losses (as explained in the next section). We next discuss metrics for measuring network design and the effectiveness of a reconnection method.

\subsection{Definition of metrics}\label{subsec:metrics}

To define metrics for reliability and self-healing capabilities in a network $G= (V,E)$, we use additional data about the demand and flow in $G$. Each vertex $v\in V$ has a demand $w(v) \in \mathbb{R}_+$. For each $e \in E$, let $p(e)$ be the probability of an event which causes edge $e$ to become inactive. Given a starting tree $T$ which connects the vertices in the network to the supply vertex (root) $r\in V$, for all $e\in T$ let $D(e)$ be the set of vertices downstream of $e$, i.e., all $v$ such that the unique path in $T$ between $r$ and $v$ includes $e$. For all $e\in T$, let the total flow $f(e) \in \mathbb{R}_+$ through $e$ be the total demand of the vertices in $D(e)$. Let the set of non-tree edges, that are initially inactive (e.g., switches), be $S = E \setminus T$. Finally, let $t(e) = \min_{s\in S: s\text{ covers }e} \sigma(s)$ be the time taken to reconnect the network and restore flow after $e$ becomes inactive, where $\sigma:S\longrightarrow [|S|]$ defines the waiting times for each non-tree edge $s\in S$. Note that $w$, $p$, $D$, and $f$ are functions of the network configuration, while $t$ depends on the network restoration methods, i.e., the choice of $\sigma$.

We define the reconnection time metric {\sc R-Time} to be the expected time between an event and the reconnection of the network, without considering the demand at each vertex: 
\[
\text{{\sc R-Time}} = \frac{\sum_{e \in T} p(e)t(e)}{\sum_{e \in T} p(e)}.
\]

We also define a metric weighted by the demand $w(v)$ of the affected vertices, which recovers the commonly used System Average Interruption Duration Index (SAIDI) reliability index for power systems~\cite{IEEE_indices}: 
\[
\text{SAIDI} = \frac{\sum_{v\in V} w(v)\sum_{e:v\in D(e)}p(e)t(e)}{\sum_{v\in V} w(v)} = \frac{\sum_{e\in T}f(e)p(e)t(e)}{\sum_{v\in V} w(v)}.
\]
For instance, if the weight $w(v)$ represents the number of consumers at vertex $v$, then SAIDI is the expected duration of time that an average consumer is without flow.

Finally, we consider a base measure of loss in the network, which is not affected by the reconnection methods but can be a substantial operational cost and is a critical factor in network design. Each edge $e\in E$ has a quantity $r(e)\in\mathbb{R}_+$, which depends on the physical properties of the corresponding object in the network (e.g., the electrical resistance of a wire). Then the total loss in the network is the quantity
\[
\text{{\sc Energy}} = \sum_{e \in T} r(e) f(e)^2.
\]
For example, in an electrical network, if we assume that vertex weight is proportional to power consumption, by Ohm's law the energy loss in any edge $e\in T$ is $r(e)f(e)^2$. (We do not consider more complex nonlinearities, e.g., AC power flow in an electric distribution network.) Similar concerns arise in gas and water networks, where when using the classical Darcy-Weisbach equation, the pressure loss for liquid flow in pipes is proportional to the length of the pipe and to the square of the velocity of flow~\cite{mueller2019comparison,eck2015quadratic}.\footnote{Gas and water networks do not always have an operational constraint of acyclicity, but local search methods could be adapted for multi-objective optimization on these.}

\subsection{The {\sc Min Reconnection Time (MRT)} problem}\label{subsec:MRT}

SAIDI and {\sc R-Time} can be viewed as special instances of the objective of a new problem, {\sc Min Reconnection Time (MRT)}. The goal of MRT is to find an ordering of the inactive edges that minimizes a weighted expected time for $T$ to be reconnected after an edge becomes inactive. That is, given a network $G=(V,E)$, a spanning tree $T$ over $V$ with edge weights $b: T \longrightarrow \mathbb{R}_+$, and a set of inactive edges $S = E\setminus T$, the MRT problem is to find an ordering $\sigma:S\longrightarrow [|S|]$ that minimizes
\begin{align}
\text{(MRT):} \quad \min_{\sigma: S \longrightarrow [|S|]} \sum_{e \in T} b(e) \min_{s \text{ covers } e} \sigma(s).
\end{align}

We set the tree edge weights $b(e)$ to reflect our choice of instance and metric. Setting $b(e)=p(e)$, the resulting objective is the expected length of time to network reconnection, while setting $b(e)=p(e)f(e)$ yields SAIDI (up to the constant factor of the denominator of those objectives). Therefore, any methods and results for solving MRT apply to the specific problems of finding a reconnection order to minimize either of these metrics.

\section{Connections to the Minimum Linear Ordering Problem}\label{sec:mssc}

In this section, we show that MRT is a special case of the Min Sum Set Cover problem, itself a special case of the Minimum Linear Ordering Problem on supermodular functions. We show that MRT is NP-hard using a reduction from the Tree Augmentation Problem. For ease and clarity of language, throughout the remainder of the paper we refer to the initial non-tree {\it self-healing} edges in $S=E\setminus T$ as ``switches,'' following the terminology for electricity distribution networks. 

\subsection{Notation and Preliminaries} 
We begin by reviewing some preliminaries needed for this section. A set function $f: 2^V\longrightarrow \mathbb{R}$ on a ground set $V$ is \textit{supermodular} if for all $A,B\subseteq V$, $f(A\cup B)+f(A\cap B)\geq f(A)+f(B)$. A \textit{hypergraph} $H= (V,E)$ is a pair of a vertex set $V$ and a hyperedge set $E$, where each hyperedge $e\in E$ is a nonempty subset (of any size) of $V$. This generalizes the notation of graphs $G=(V,E)$ where each edge $e \in E$ is a subset of $V$ containing exactly two vertices, i.e., $e = (u,v)$ for $u,v\in V$. A hypergraph is said to be \textit{$c$-uniform} if each hyperedge contains exactly $c$ vertices. We denote linear orderings of a finite set $V$ by $\sigma: V \longrightarrow [|V|]$, where $\sigma(i)$ is the position of element $i \in V$.

We show that MRT is an instance of the {\sc Minimum Linear Ordering Problem} (MLOP) on supermodular functions, in particular the {\sc Min Sum Set Cover} (MSSC) problem introduced by~\cite{Feige2002}. We first define the general form of MLOP:

\begin{definition}
Given a ground set $S$ and a function $f: 2^S\longrightarrow \mathbb{R}_+$, the MLOP is to find a linear ordering $\sigma$ of $S$ that minimizes $\sum_{i=1}^{|S|} f(\{s\in S: \sigma(s)\leq i\})$.
\end{definition}

To see how MRT can be viewed as a special case of MLOP, consider a spanning tree $T$ of the underlying network $G = (V, E)$. The goal is to order the ground set of switches $S=E\setminus T$ to minimize the time for edges in $T$ to be covered. For each switch $s$, let $U_s$ be the set of tree edges covered by $s$. Given a set of switches $W\subseteq S$, let $f(W)$ be the cardinality of the set of tree edges not covered by $W$, i.e., $f(W) = |T \setminus \cup_{s\in W} U_s|$. This function can be shown to be supermodular. At each time $t=i$, any edge not covered by the first $i-1$ switches adds 1 unit to the objective function. As in MRT, the value of the set function $f$ can be weighted by $b(e)$, thus recovering {\sc R-Time} or SAIDI as the objective. Any ordering that minimizes the MLOP objective also minimizes MRT.\footnote{There is an additive factor of $|S|$ in the objective that is not accounted for in the MLOP instance but does not change the optimal solution.}

Another special case of the MLOP is the Min Sum Set Cover (MSSC) problem, in which given a set of subsets $S_1, \hdots, S_k \subseteq S$, we find an ordering of the $S_i$ such that the total waiting time for each element of the ground set $S$ to be covered is minimized. As an instance of the MLOP, this problem uses the function $f(U) = |S \setminus \cup_{i\in U} S_i|$ for $U \subseteq V$, which again is supermodular. We may view this instance using the hyperedges of a hypergraph $H= (V,E)$, where each hyperedge $e\in E$ is a nonempty subset (of any size) of the vertex set $V$. Taking the ground set to be the edge set $E$, the MSSC problem is formally stated as:
\begin{definition}
The {\sc MSSC} problem, given a hypergraph $H=(V,E)$, is to find a linear ordering $\sigma$ of the vertices $V$ so that the total waiting time for each edge to be covered is minimized, i.e., $\sum_e \min_{v\in e} \sigma(v)$ is minimized.
\end{definition}

Again, note that one can write an instance of {\sc MRT} as an instance of {\sc MSSC} by taking $H$ to be the hypergraph with a vertex for each switch $s\in S=E\setminus T$ and a hyperedge $\{s\in S:T-e+f \text{ is connected}\}$ for each tree edge $e\in T$. However, MRT instances form a strict subset of MSSC instances, so showing MRT is NP-hard requires additional argument.

\begin{lemma}\label{lem:MSSCspecialcase}
There exist instances of {\sc MSSC} which are not instances of {\sc MRT}.
\end{lemma}
\emph{Proof.}
Let $H=(\{1,2,3,4\}, \{\{1,2,3\},\{1,4\},\{2,4\},\{3,4\}\})$. Then in any corresponding instance of {\sc MRT}, for a set of switches $S$ corresponding to the hyperedges of $H$ to exist, tree edges $e_1,e_2,e_3$ must lie on a path, with a fourth tree edge $e_4$ incident to each of them. However, this is impossible. 
\QED 

NP-hardness also may not hold in special cases of MRT, for instance if $G$ is planar.\footnote{Analogously, even though finding the maximum weighted cut in a graph is NP-hard~\cite{goemans1995improved}, it is solvable in polynomial time for planar instances~\cite{dorfman1972finding}.} There may also be more effective heuristics or stronger approximation guarantees if $G$ is planar~\cite{whited1990reliability}.

\subsection{NP-hardness of Min Reconnection Time}
To show the NP-hardness of the Min Reconnection Time problem (MRT), we show a reduction from another NP-hard problem called the {\sc Tree Augmentation} problem, or TAP~\cite{Cohen2013}. Given $G=(V,E)$ and a spanning tree $T\subseteq E$ on $V$, the objective of TAP is to find a minimum-cardinality set of edges $F \subseteq E$ such that $T\cup F$ is 2-edge-connected. Compare this to MRT, in which an ordering $\sigma$ (a solution to MRT) is guaranteed to restore power by time $k$ if the subgraph of $G$ containing the edges of $T \cup \{s\in S: \sigma(s)\leq k\}$  is a 2-edge-connected graph. Thus MRT can be viewed as a timed version of TAP, much as MSSC is the timed version of ordinary {\sc Set Cover}: the latter seeks to find the minimum \textit{number} of subsets to cover a given set of elements~\cite{karp1972reducibility}.

\begin{theorem}\label{thm:hardness} The tree augmentation problem can be reduced in polynomial time to the {\sc Min Reconnection Time} problem, and the latter problem is therefore NP-hard.
\end{theorem}
\emph{Proof.}
Let $G=(V,E)$ with spanning tree $T$ be an instance of TAP, and let $n=|V|, m=|E|$.

We modify the instance to create a new graph $G'= (V',E')$ with a new spanning tree $T'$. Add the $2mn$ vertices $v_{i,j}$, where $i\in [mn],j\in [2]$, setting $V' = V\cup\{v_{1,1},v_{1,2}, \dots, v_{mn,1},v_{mn,2}\}$. Select an arbitrary vertex $v\in V$, and let $P_T=\{\{v,v_{1,1}\}, \{v_{1,1},v_{1,2}\}, \dots, \{v,v_{mn,1}\}, \{v_{mn,1},v_{mn,2}\}\}$ be a set of edges which constitute $mn$ disjoint 2-paths. Let $P_S = \{\{v,v_{1,2}\}, \dots, \{v,v_{mn,2}\}\}$ be a set of $mn$ edges, and note that for $1\leq i\leq mn$, $\{v,v_{i,2}\}$ covers the edges $\{v,v_{1,1}\}, \{v_{1,1},v_{1,2}\}$. Finally, let $E'= E\cup P_T\cup P_S$, let $T' = T\cup P_T$, and give all edges $e\in T$ weight $b(e)=1$ and all edges $e\in P_T$ weight $b(e)=1/2$.

Let $\sigma$ be an optimum solution to {\sc MRT} on $G'=(V',E')$ with tree $T'$. The newly added switches $\{v,v_{1,2}\},\dots,\{v,v_{mn,2}\}$ are the only switches that cover edges in $P_T$. Furthermore, these switches each cover tree edges of total weight 1, and since all edges in $T$ have weight 1, without loss of generality we may assume that the switches in $P_S$ are the last $mn$ switches in $\sigma$ that contribute a positive quantity to the objective $OBJ(\sigma)$ of $\sigma$ as a solution to {\sc MRT}. Let the switches chosen before those $mn$ switches be $s_1,\dots,s_k = \sigma^{-1}(1),\dots,\sigma^{-1}(k)$, and for $1\leq i \leq k$, let $S_i=\{e\in T: s_i\text{ covers } e\}$. Then the objective value of $\sigma$ is
\begin{align*}
OBJ(\sigma) &= \sum_{i=1}^{k} i \left| S_i\setminus \bigcup_{j=1}^{i-1}S_j\right| + \sum_{\ell=k+1}^{k+mn} 1\\
 &= \sum_{i=1}^{k} i \left| S_i\setminus \bigcup_{j=1}^{i-1}S_j\right| + mn\cdot k + \frac{mn(mn+1)}{2}.
\end{align*}
Suppose there is a solution $F=\{f_1,\dots,f_{k^*}\}\subseteq S$ to TAP on the original instance $G,\ T$ such that $k^*<k$. We construct another solution $\sigma^*$ to MRT on $G'$ and $T'$ as follows.

For $1\leq i \leq k^*$, let $F_i=\{e\in T: f_i\text{ covers } e\}$. Without loss of generality, assume that $\left| F_i\setminus \bigcup_{j=1}^{i-1}F_j\right|\geq \left| F_{i+1}\setminus \bigcup_{j=1}^{i}F_j\right|$ for all $i\in [k^*-1]$ (if not, we can permute the indices of the $f_i$ to make this so). Let $\sigma^*$ be a permutation of $E'\setminus T'$ such that
\[
\sigma^{*-1}(i) =
\begin{cases}
f_i, \text{ if } i\leq k^* \\
\{v,v_{i-k^*,2}\}, \text{ if } k^*+1\leq i\leq k^*+mn.
\end{cases}
\]
Then since $\sum_{i=1}^{k^*}\left| F_i\setminus \bigcup_{j=1}^{i-1}F_j\right|=\left| \bigcup_{i=1}^{k^*}F_j\right|=|T|=n-1$, we have that
\[
\sum_{i=1}^{k^*} i\left| F_i\setminus \bigcup_{j=1}^{i-1}F_j\right| \leq \sum_{i=1}^{k^*} i\frac{n-1}{k^*} = \frac{k^*(k^*+1)}{2}\frac{n-1}{k^*} =\frac{(k^*+1)(n-1)}{2} \leq mn.
\]
Therefore, 
\begin{align*}
    OBJ(\sigma^*) = \sum_{i=1}^{k^*} i\left| F_i\setminus \bigcup_{j=1}^{i-1}F_j\right| + \sum_{\ell=k^*+1}^{k^*+mn} 1 &\leq mn + mn\cdot k^* + \frac{mn(mn+1)}{2} \\
    &\leq mn\cdot k + \frac{mn(mn+1)}{2} < OBJ(\sigma),
\end{align*}
which would contradict the optimality of $\sigma$. Therefore, there is no such solution $\{f_1,\dots,f_{k^*}\}$ to TAP on $G$ and $T$ using $k^*<k$ edges, and so $\{s_1,\dots,s_k\}$ is an optimum solution to the original TAP instance. Since $G'$ and $T'$ have size polynomial in the size of $G$ and $T$, this gives a polynomial-time reduction of TAP to {\sc MRT}, and thus {\sc MRT} is NP-hard.
\QED

\section{Approximation algorithms for MRT}\label{sec:approx} 
We next discuss methods for solving MSSC and MRT with provable approximation guarantees. We first show that the greedy algorithm is not optimal for MRT, and then develop new fractional kernels for extending the alpha-point rounding method for MSSC introduced in~\cite{Bansal2021}. 
\subsection{The greedy algorithm}
Since an instance of MRT is an instance of MSSC, a natural question is to ask if the greedy algorithm \cite{Feige2002} gives a good approximation factor for the MRT as well. For each $i$ from 1 to $|S|$, the greedy algorithm chooses 
\[
\sigma^{-1}(i) = \argmax_{s\in S}\sum_{\substack{e \in T\\ s\text{ covers }e\\\sigma_j\text{ does not cover }e, \forall j<i}} b(e);
\]
that is, $\sigma^{-1}(i)$ is the edge in $S$ which covers the uncovered tree edges with greatest total weight. \cite{Feige2002} showed that the greedy algorithm is a 4-approximation for general MSSC, and that this bound is tight.

There exist instances of MRT where the greedy algorithm does not find an optimal solution, and we construct instances with a gap of 4/3. We suspect that the greedy algorithm achieves an approximation ratio of strictly less than 4 on MRT instances.  

\begin{figure}[t]
\begin{center}
\begin{tikzpicture}[scale=1,every node/.style={draw=black,circle}]
    \node (s) at (0,0) {s};
    \node (1) at (-1,-2) {$v_1$};
    \node (2) at (-2,-4) {$v_2$};
    \node (3) at (1,-2) {$v_3$};
    \node (4) at (2,-4) {$v_4$};
    \draw[->] (s) to (1);
    \draw[->] (1) to  (2) ;
    \draw[->] (s) to (3);
    \draw[->] (3) to (4);
  \tikzset{switch/.style={black,fill=none,draw=none,midway}};
    \draw[dashed] (s) to[out = 200, in = 110] node[switch,left] {$s_2$} (2);
\draw[dashed] (s) to[out = -20, in =70] node[switch,right] {$s_3$} (4);
    \draw[dashed] (1) to node[switch,below] {$s_1$} (3);
\end{tikzpicture}
\includegraphics[width=.33\linewidth]{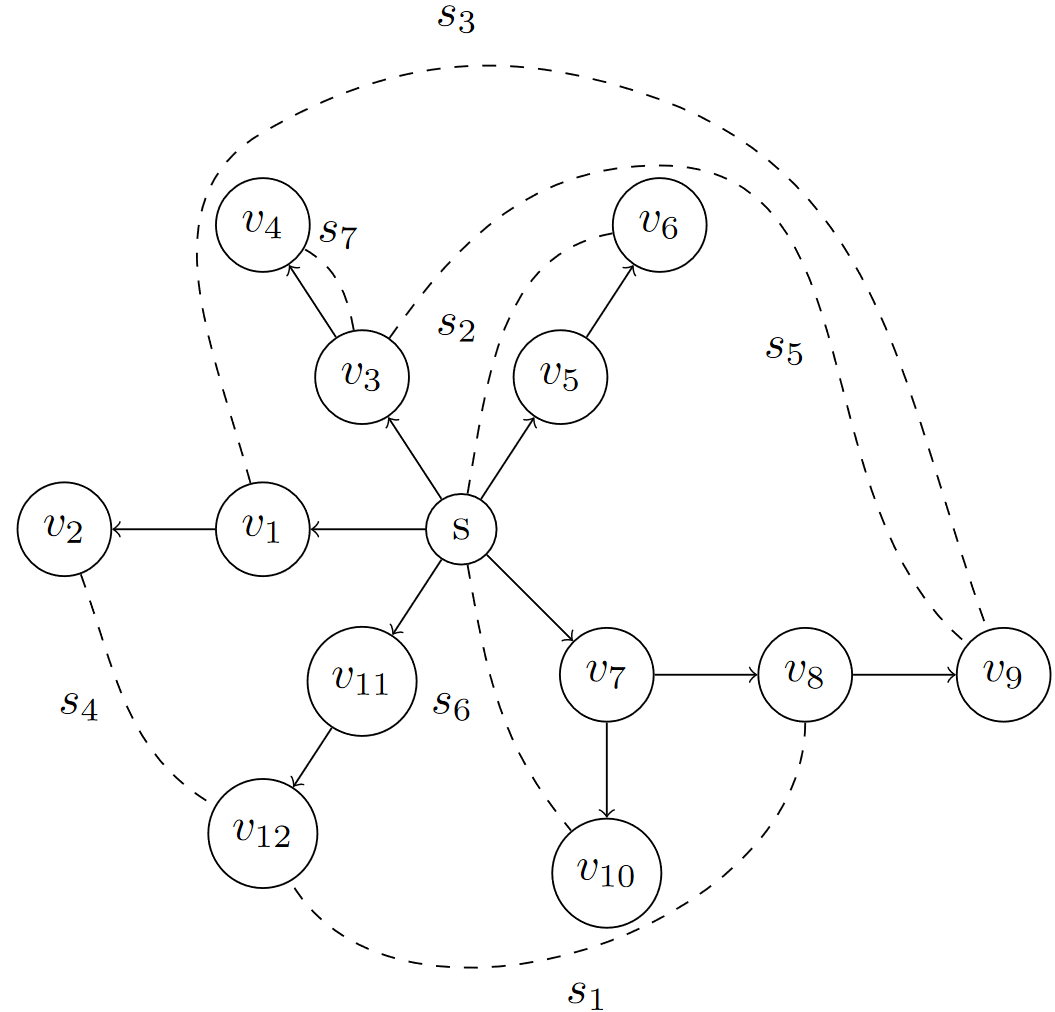}
    \end{center}
    \caption{Two networks where the greedy algorithm for MRT does not achieve the optimum value.}
    \label{fig:greedyexample}
\end{figure}
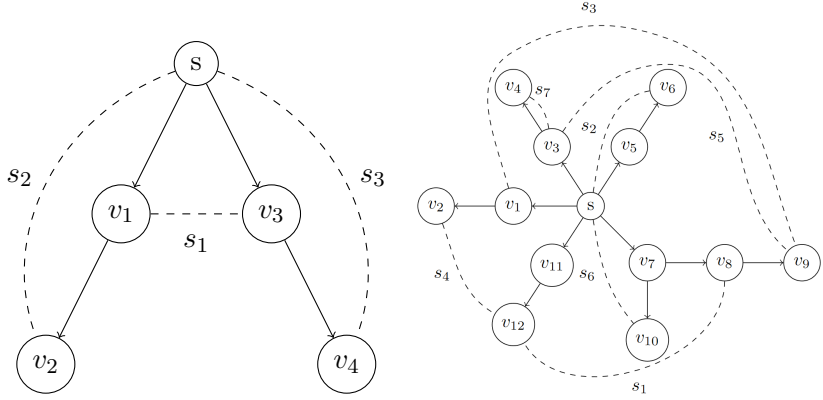

\vspace{0.2cm}
\noindent 
\textbf{Example 1:} Consider the leftmost network in Figure~\ref{fig:greedyexample}. Edges in $T$ are shown solid, while edges in $S$ are dashed. Let $\varepsilon=1$ and let $p(e)=1$ for all $e\in T$. Then an optimal solution to MRT on this is to take the edges of $S$ in order $\sigma = (s_2,s_3,s_1)$, yielding objective value $3/2$. However, the greedy algorithm, breaking ties lexicographically, will take $\sigma = (s_1,s_2,s_3)$, yielding objective value $7/4$, a multiplicative gap of $7/6$.

\vspace{0.2cm}
\noindent 
\textbf{Example 2:} Another instance when the greedy algorithm is not optimal is depicted in Figure~\ref{fig:greedyexample} (right). Again, let $\varepsilon=1$ and let $p(e)=1$ for all $e\in T$. Then an optimal solution to MRT is $\sigma = (s_4,s_5,s_2,s_6,s_7,s_1,s_3)$, yielding objective value $27/12$. However, the greedy algorithm, breaking ties lexicographically, selects $\sigma = (s_1,s_2,s_3,s_4,s_5,s_6,s_7)$, yielding an objective value of $36/12$, a 4/3 multiplicative gap to the optimum. In this network, the edge $e$ between $v_3$ and $v_4$ is covered only by edge $s_7$, so any optimal ordering on $S$ can take $s_7$ as the last edge $\sigma^{-1}(i)$ that covers an edge not covered by $\sigma^{-1}(1),\dots,\sigma^{-1}(i-1)$. In this sense, the existence of $e$ does not make solving the MRT problem any more difficult. However, if $e$ is removed from $G$, the multiplicative gap shrinks to 29/22.

\subsection{Kernel-based rounding for 
\texorpdfstring{$c$}{c}
-uniform MRT instances}
We next show how to improve the approximation ratio for Min Sum Set Cover, including MRT, by extending a recent technique called kernel-based alpha-point rounding ~\cite{Bansal2021}. Unlike the approach taken by Bansal et. al, we use {\it fractional kernels} within a kernel-based randomized rounding approach~\cite{phillips1998minimizing,hall1997scheduling}, which may be of independent interest. In particular, we improve approximations for $c$-uniform instances of MSSC.

Given an MRT instance with every switch covering at most $c$ edges of the underlying operational radial network, we first convert such an instance to a minimum sum set cover (MSSC) instance on a $c$-regular hypergraph (by analogy to the 3-layered zones of protection \cite{blackburn2015,glover2017}, we expect that tree edges are usually covered by a small number of switches). For these instances, our new kernels interpolate the above mentioned approximation factors from $16/9 \approx 1.778$ to 4, dependent on $c$. Therefore, these kernels strictly improve the state-of-the-art factors for MSSC instances. 

The kernel-based $\alpha$-point rounding method begins by finding a solution $x$ to an LP relaxation of the following integer program.\footnote{Our formulation differs slightly from that of Bansal et al., specifically in the covering constraints~\eqref{eq:covering1} and ~\eqref{eq:covering2}. We use this modification for its improved computational performance. Otherwise the modified formulation does not affect the algorithm, as our LP is polynomial size and uses the same decision variables.} The formulation encodes MSSC using binary variables $x_{v,t}$ to denote whether vertex $v$ is chosen at time step $t$. 
\begin{align}
    \mathrm{minimize}\qquad & \sum_{t,e}t\cdot b(e)\cdot u_{e,t}\\
     \mathrm{subject\ to}\qquad & \sum_{v\in e} x_{v,t}-u_{e,t} \geq 0\qquad \forall t \text{ time periods, } e \text{ edges},\\
    \qquad & \sum_{t} u_{e,t} =1\qquad \forall e \text{  edges}, \label{eq:covering1}\\
    \qquad & \sum_{v} x_{v,t} =1\qquad \forall t \text{ time periods}, \label{eq:covering2}\\
    &\text{all  $u_{e,t},x_{v,t}\geq 0$ for vertices $v$, time $t$ and edges $e$.}
\end{align}
If the variables $u_{e,t},x_{v,t}$ are binaries, then indices $t\in [|V|]$ denote the index in the selection order $\sigma$, so that $\sigma(v) = t$ if and only if $x_{v,t} = 1$. Variables $u_{e,t}$ denote the time $t$ at which edge $e$ is covered. Each edge $e$ must be covered at some time $t$, and at each time $t$, at most one vertex can be chosen.

The LP relaxation can be solved in polynomial time. We apply a kernel $K(t,t')$ to $x$ to obtain new weights $z$. Applying $\alpha$-point randomized rounding, we obtain a provisional schedule of selection times $\tau$ for all vertices $v$, then convert this list to our desired ordering $\sigma$ by taking the vertices in the order in which they appear in $\tau$, breaking ties at random. Pseudocode for this algorithm appears in Algorithm~\ref{alg:rrmrt}. 

\cite{Bansal2021} obtained a 4-approximation for MSSC using $K = 2/t \cdot \mathds{1}[t\geq t']$, and a 16/9-approximation for the minimum sum vertex cover problem (MSVC) using $K = 4 \frac{t'(t'+1)}{t(t+1)(t+2)} \cdot \mathds{1}[t\geq t']$. In contrast, to attain tight guarantees for MSSC on $c$-uniform hypergraphs, we consider a new kernel of the form $K(t,t')=\frac{2c}{c+1} \frac{t'^{2/(c-1)}}{\sum_{i=1}^{t} i^{2/(c-1)}} \cdot \mathds{1}[t\geq t']$,
which is a discrete analog of negative power functions. This kernel spreads the weight of relaxed LP solutions more aggressively, which is necessary to ensure that small fractional contributions in weight (which are more frequent for small $c$) contribute to a good tentative schedule. The analysis of the gap between the LP objective and the cost of the tentative schedule uses discrete calculus-like manipulations of the kernels, and requires close bounds on their sums. These are significantly more intricate when dealing with the fractional powers of $t$ used in our kernels. Our bounds on the total cumulative weight assigned to the vertices of an edge $e$ require additional convex analysis, including two more inequalities  (Lemmas~\ref{lem:isum} and~\ref{lem:HHlike}).



\begin{algorithm}[!t]\footnotesize
 \caption{\textsc{Kernel $\alpha$-point Rounding for MSSC \cite{Bansal2021}}}\label{alg:rrmrt}
 \begin{algorithmic}[1]
\State \textbf{Input:} Hypergraph $H=(V,E)$; kernel $K(t,t')$.
\State \textbf{Output:} An ordering $\sigma:V\longrightarrow [|V|]$ of $V$.
\State $x\leftarrow$ Optimal fractional solution for MRT.
\State $z_{s,t}\leftarrow \sum_{t'}K(t,t')x_{v,t'}\ \forall v\in V,1\leq t \leq |V|$
\For{$v\in V$}
    \State Sample $\alpha_v\sim [0,1].$
    \State $\tau_{v} \leftarrow \text{the earliest time $t$ for which }\sum_{t'\leq t} z_{v,t'} \geq \alpha_v$
\EndFor
\State \Return an ordering $\sigma$ of $V$, scheduling vertices according to $\tau$ and breaking ties at random.
\end{algorithmic}
\label{alg:mrt_rr}
\end{algorithm}

\begin{theorem}\label{thm:alphaallc}
Let $H=(V,E)$ be a $c$-uniform hypergraph, where $c\geq 3$. Then applying $\alpha$-point rounding with kernel 
$K(t,t')=\frac{2c}{c+1} \frac{t'^{2/(c-1)}}{\sum_{i=1}^{t} i^{2/(c-1)}} \cdot \mathds{1}[t\geq t']$
is a polynomial-time $(\frac{2c}{c+1})^2$ approximation for the Min Sum Set Cover problem on $H$.
\end{theorem}

\emph{Proof.}
As defined in the algorithm, let $z_{s, t} = \sum_{t^\prime} K(t,t^\prime) x_{s,t^\prime}$ be the weight assigned to any vertex $s$ and time period $t$ after rounding. Fix a hyperedge $e\in E$. Let $e=\{v_1,\dots,v_c\}$. Let $z_{e,t}=\sum_{i=1}^{c} z_{v_i,t}$ be the total weight of vertices covering $e$ at exactly time $t$, after applying the kernel $K$ to the LP solution $x$, and let $z_{e,<t} = \sum_{t'<t}z_{e,t'}$. Let $q_t(e) = \PR[e \text{ is not scheduled before time $t$ in $\tau$}]$, and define $p_t(e) = ((1-z_{e,<t}/c)^c)_+$, which we will show is an upper bound on the former probability.

Define the cost of covering $e \in E$ in the LP ($c_x(e)$), a bound on the time $e$ is scheduled in the tentative order $\tau$ ($c_z(e)$), and the cost of covering $e$ in the final order $\sigma$ ($c_{\sigma}(e)$) as: 
\begin{enumerate}
\item[(i)] LP cost: $c_x(e)=\sum_t (1-\sum_{t'<t}\sum_{i=1}^{c}x_{v_i,t'})_+$, where $(\cdot)_+ = \max \{0,\cdot\}$.
\item[(ii)] Upper bound on expected schedule time in $\tau$: $c_{z}(e) = \sum_{t}p_t(e)$.
\item[(iii)] Cost in $\sigma$: $c_{\sigma}(e) = \E[\min\{t:\sigma^{-1}(t) \in \{v_1,\dots,v_c\}\}]$. 
\end{enumerate} 

The LP objective is $\sum_{e \in E}c_x(e)$ and the objective value of the ordering returned by $\alpha$-point rounding is $\sum_{e\in T}c_{\sigma}(e)$. Therefore, it suffices to show that $(\frac{2c}{c+1})^2c_x(e) \geq c_{\sigma}(e)$ for all $e\in E$. To do so, we separately show that $(\frac{2c}{c+1}) c_x(e) \geq   c_z(e)$ and $(\frac{2c}{c+1})c_z(e) \geq   c_{\sigma}(e)$.

\noindent 
{\bf Part A. We would like to show first that $(\frac{2c}{c+1}) c_z(e) \geq c_{\sigma}(e)$.} 
The following argument can be generalized to a variety of kernels for $\alpha$-point rounding.\footnote{In particular, claim 1 requires that for some $\beta\geq 1$, $\sum_{t'}K(t,t')\leq \beta$ for all $t$; claim 2 requires that for all $t'\geq 1$, $\sum_{t}K(t,t')\geq 1$; and claims 3, 4, and 5 do not otherwise depend on $K$.} 
We apply it with $\beta=\frac{2c}{c+1}.$ 

\noindent{\bf Claim 1. For the kernel $K(t^\prime, t) = \frac{2c}{c+1} \frac{t'^{2/(c-1)}}{\sum_{i=1}^{t} i^{2/(c-1)}} \cdot \mathds{1}[t\geq t']$, $\sum_{v} z_{v,t} \leq 2c/(c+1)=\beta$.}
    
    We have that for any $t$,
\begin{align*}
\sum_{v} z_{v,t} &= \sum_{v}\sum_{t'}K(t,t')x_{v,t} = \sum_{t'}K(t,t')\sum_v x_{v,t} \\
&\leq \sum_{t'}K(t,t') = \frac{2c}{c+1}\sum_{t^\prime=1}^t\frac{t'^{2/(c-1)}}{\sum_{i=1}^{t} i^{2/(c-1)}} \\
&=\frac{2c}{c+1}\frac{\sum_{t^\prime=1}^t t'^{2/(c-1)}}{\sum_{i=1}^{t} i^{2/(c-1)}} = \frac{2c}{c+1}.
\end{align*}

\noindent{\bf Claim 2. The quantity $p_t(e)=((1-z_{e,<t}/c)^c)_+$ upper-bounds $q_t(e)$ and is 0 for sufficiently large $t$:}
For all $v$, we have that
    \begin{align*}
    \sum_t z_{v,t} &= \sum_t \sum_{t'} K(t,t') x_{v,t'} = \sum_{t'} x_{v,t'} \sum_t K(t,t')\\
    &=\sum_{t'} x_{v,t'} \frac{2c}{c+1}t'^{2/(c-1)}\sum_{t\geq t'} \frac{1}{\sum_{i=1}^{t} i^{2/(c-1)}} \\
     &>\sum_{t'} x_{v,t'} \frac{2c}{c+1}t'^{2/(c-1)}\sum_{t\geq t'} \frac{1}{\int_{0}^{t+1} i^{2/(c-1)} di} \\
 &=\sum_{t'} x_{v,t'} \frac{2c}{c+1}t'^{2/(c-1)}\sum_{t\geq t'} \frac{c+1}{c-1}(t+1)^{-(c+1)/(c-1)} \\
  &>\sum_{t'} x_{v,t'} \frac{2c}{c-1}t'^{2/(c-1)}\int_{t'+1}^\infty (t+1)^{-(c+1)/(c-1)} dt \\
    &= \sum_{t'} x_{v,t'}\cdot c \left(\frac{t'}{t'+1}\right)^{2/(c-1)}\\
    &\geq \sum_{t'} x_{v,t'} \cdot 3 (1/2)^{1}> \sum_{t'} x_{v,t'} \cdot 1 = 1, \numberthis\label{eq:claim2}
    \end{align*}
    where (\ref{eq:claim2}) follows since without loss of generality we may assume that in the solution of the linear program for MSSC, each vertex $v$ is selected with total weight 1 across all time steps.
    
    For all $v$, since $\sum_{t=1}^{\infty} z_{v,t}>1$, for sufficiently large $t$ we have that $z_{v,<t}=\sum_{i=1}^{t-1} z_{v,i}\geq 1$. Therefore, $z_{e,<t} = \sum_{i=1}^{c} z_{v_i,<t}\geq c$ for sufficiently large $t$, so $p_t(e) = 0$ for sufficiently large $t$. 
    
    Moreover, by the arithmetic mean-geometric mean inequality, we have that $q_t(e)$ is bounded above by $p_t(e)$, since
\begin{align*}
q_t(e)&= \prod_{i=1}^{c}(1-z_{v_i,<t}) \leq \left(1-\sum_{i=1}^{c} z_{v_i,<t}/c\right)^c_+ =p_t(e).
\end{align*} 

\noindent{\bf Claim 3. Upper-bound the expected cost of $e$ in the final ordering $\sigma$:}  
We claim
    \begin{equation}\label{eqn:sigmabound}
        \E[c_\sigma(e)] \leq \sum_t (q_t(e) - q_{t+1}(e))(1+\beta(t-1)-z_{e,<t}+(\beta-z_{e,t})/2).
    \end{equation}
    Each term of this sum corresponds to the event that $e$ is scheduled at exactly time $t$ in $\tau$. This happens with probability $q_t(e) - q_{t+1}(e)$. The quantity $S(t)=1+\beta(t-1)-z_{e,<t}+(\beta-z_{e,t})/2$ is an upper bound on the expected number of vertices selected in $\sigma$ before and including the one that covers $e$. By Claim 1, in this quantity:
    \begin{itemize}
        \item 1 corresponds to the first vertex which covers $e$.
        \item $\beta(t-1)-z_{e,<t}$ is an upper bound on the expected number of vertices scheduled in $\tau$ before $e$ is covered. Indeed, that quantity is $\sum_{v\not\in e}z_{v,<t} \leq \beta(t-1) - z_{e,<t}.$
        \item In expectation, at most $\beta$ vertices are scheduled at exactly time $t$ in $\tau$, of which $(\beta-z_{e,t})$ do not cover $e$. In expectation, at most half of these vertices will be scheduled in $\sigma$ before $e$ is covered since ties are broken at random, yielding the bound of $(\beta-z_{e,t})/2$.
    \end{itemize}

$S(t)$ is non-decreasing with respect to $t$. For each $t$, replacing $q_t(e)$ with $p_t(e)$ increases the right-hand side of Equation~\eqref{eqn:sigmabound} by $(p(t)-q(t))(S(t+1)-S(t))>0$. Therefore, since $p_t(e)$ is positive for a finite number of $t$ (Claim 2), replacing all $q_t(e)$ with the corresponding $p_t(e)$ can only increase the sum, and therefore
\[
      \E[c_\sigma(e)] \leq \sum_t (p_t(e) - p_{t+1}(e))(1+\beta(t-1)-z_{e,<t}+(\beta-z_{e,t})/2).
\]

Note that $\E[c_z(e)] = \sum_t p_t(e) =\sum_t t(p_t(e)-p_{t+1}(e))$, so $\E[c_\sigma(e)]\leq \beta\E[c_z(e)]$ holds if and only if each of the following inequalities holds. We prove the last of these in Claim 4.
\begin{align*}
    \sum_t (p_t(e) - p_{t+1}(e))(1+\beta(t-1)-z_{e,<t}&+(\beta-z_{e,t})/2)\leq \beta \sum_t t(p_t(e)-p_{t+1}(e)) \\
    \sum_t (p_t(e) - p_{t+1}(e))(1 -z_{e,<t}-z_{e,t}/2) &\leq  \beta/2 \sum_t (p_t(e) - p_{t+1}(e)) \\ 
    (1-\beta/2) \sum_t (p_t(e) - p_{t+1}(e)) &\leq \sum_t (p_t(e) - p_{t+1}(e))(z_{e,<t}+z_{e,t}/2)\\
    1-\beta/2 &\leq \sum_t (p_t(e) - p_{t+1}(e))(z_{e,<t}+z_{e,t}/2)\\
    1/2 &\leq \sum_t (z_{e,t}/2) (p_t(e)+p_{t+1}(e))\\
    1 &\leq \sum_t z_{e,t}(p_t(e) + p_{t+1}(e)).
\end{align*}

\noindent{\bf Claim 4. Lower-bound $\sum_t z_{e,t}(p_t(e) + p_{t+1}(e))$:} Note that $p_t(e)$ is a convex function of $z_{e,<t}$, as it is 0 when $z_{e,<t}\geq c$. Therefore, for all $a<b$ we have that $(b-a)(p_t(a)+p_t(b))\geq 2\int_a^b p_t(u)\ du$. In particular, for all $t$ we have that 
\[
z_{e,t}(p_t(e)+p_{t+1}(e)) = (z_{e,t+1}-z_{e,t})(p_t(e)+p_{t+1}(e)) \geq 2\int_{z_{e,<t}}^{z_{e,<t+1}} p_u(e)\ du,
\]
and so
\begin{align*}
\sum_t z_{e,t}(p_t(e) + p_{t+1}(e)) &\geq 2\sum_t \int_{z_{e,<t}}^{z_{e,<t+1}} p_u(e)\ du 
\\
&= 2 \int_{0}^{c} (1-u/c)^c\ du = 2c/(c+1) >1.
\end{align*}

\noindent{\bf Claim 5. Finally, we claim that $(2c/(c+1)) c_z(e) \geq c_{\sigma}(e)$:} Claim 3 holds using $\beta = 2c/(c+1)$ (Claim 1). Applying Claim 4, we have that $(2c/(c+1)) c_z(e) \geq c_{\sigma}(e)$ as desired.

\ 

\noindent\textbf{Part B. We would now like to show that $(2c/(c+1)) c_x(e) \geq c_z(e)$.} 

For all vertices $v$, we have that
\[
z_{v,<t} = \sum_{t'<t}z_{v,t'} = \sum_{t'<t}\sum_{t''<t'}K(t',t'')x_{v,t''} = \sum_{t''<t} x_{v,t''} \sum_{t'=t''}^{t-1}K(t',t'').
\]

To analyze this sum, we use the following lemma, which we prove in Appendix~\ref{app:isumlemmaproof}.
\begin{lemma}\label{lem:isum} For $t\in\mathbb{N}$ and $0<p\leq 1$, $\sum_{i=1}^{t} i^p \leq\frac{p}{p+1}\frac{t^p (t+1)^p}{(t+1)^p-t^p}$.
\end{lemma}

Examining the inner sum $\sum_{q=t''}^{t-1}K(q,t'')$ more closely and applying Lemma~\ref{lem:isum} with $p=2/(c-1)$, we have that 
\begin{align*}
\sum_{q'=t''}^{t-1}K(q,t'') &= \frac{2c}{c+1}t''^{2/(c-1)} \sum_{q'=t''}^{t-1} \frac{1}{\sum_{i=1}^{q}i^{2/(c-1)}}\\
&\geq \frac{2c}{c+1}t''^{2/(c-1)}\sum_{q'=t''}^{t-1}\frac{c+1}{2}\left(\frac{1}{q^{2/(c-1)}}-\frac{1}{(q+1)^{2/(c-1)}}\right)\\
&=\frac{2c}{c+1}t''^{2/(c-1)}\frac{c+1}{2}\left(\frac{1}{t''^{2/(c-1)}}-\frac{1}{t^{2/(c-1)}}\right)\\
&=c\left(1-\frac{t''^{2/(c-1)}}{t^{2/(c-1)}}\right).
\end{align*}
Thus 
\[
z_{v,<t} = \sum_{t''<t} x_{v,t''} \sum_{t'=t''}^{t-1}K(t',t'')\geq \sum_{t''<t} c\left(1-\frac{t''^{2/(c-1)}}{t^{2/(c-1)}}\right) x_{v,t''},
\]
and so
\begin{align*}
c_{z}(e) &= \sum_{t} p_{t}(e)=\sum_t \left(1-\sum_{i=1}^{c}z_{v_i,<t}/c\right)_+^c\\
&\leq \sum_{t} \left(1-\sum_{t'<t}\left( 1-\frac{t'^{2/(c-1)}}{t^{2/(c-1)}}\right)\sum_{i=1}^{c}x_{v_i,t'}\right)_+^c.
\end{align*}

We now modify the LP solution $x$ in order to bound $c_{z}(e)$. These modifications do not decrease our lower bound  $c_{z}(e)$ and do not increase the LP cost $c_x(e)$. First, decrease $x$ as necessary so that $\sum_{i=1}^{c}x_{v_i,<t}=\sum_{t'<t}\sum_{i=1}^{c}x_{v_i,t'}\leq 1$ for all $t$. Second, for all $t$, replace each $x_{v_i,t}$ by their average value $a_t = \sum_{i=1}^{c}x_{v_i,t}/c$. As $c_z(e)$ depends only on the total $\sum_{i=1}^{c}x_{v_i,t}$ and not the distribution of weight between these vertices, this does not change $c_z(e)$. Then our goal is to show that for all nonnegative $a$ with $||a||_1=1$,
\begin{align*}
c_{z}(e) &\leq\sum_{t} \left(1-\sum_{t'<t}\left( 1-\frac{t'^{2/(c-1)}}{t^{2/(c-1)}}\right)\sum_{i=1}^{c}x_{v_i,t'}\right)_+^c \\
&\leq \frac{2c}{c+1} \sum_{t}\left(1-\sum_{t'<t}a_{t'}\right)^c = \frac{2c}{c+1}c_x(e).
\end{align*}
For each $t$, the corresponding summand $(1-d\cdot a)^c$ is positive and a convex function of $a$, since the vector $d$ given by $(d)_{t'} = t'/t \cdot \mathds{1}[t'<t]$ satisfies $||d||_\infty\leq 1$. Therefore $c_z(e)$ is a convex function of $a$. Moreover, $c_x(e)$ is a linear function of $a$. Therefore, by Fact 8 of~\cite{Bansal2021}, the quotient $c_z(e)/c_x(e)$ is maximized at an extreme point $a$ of the positive simplex, i.e.\ for some $u$, $a_u=1$ and $a_{u'}=0$ for $u'\neq u$. At this point $a$, we have that $2c/(c+1) c_x(e)=2c/(c+1) u$ and
\begin{align*}
c_z(e) &= \sum_t \left(1-\sum_{t'<t}\left(1-\frac{t'^{2/(c-1)}}{t^{2/(c-1)}}\right)a_{t'}\right)^c \\
&=  \sum_{t=1}^{u} \left(1-\sum_{t'<t}\left(1-\frac{t'^{2/(c-1)}}{t^{2/(c-1)}}\right)a_{t'}\right)^c + \sum_{t>u} \left(1-\sum_{t'<t}\left(1-\frac{t'^{2/(c-1)}}{t^{2/(c-1)}}\right)a_{t'}\right)^c \\
&= u+\sum_{t>u} \left(1-\left(1-\frac{u^{2/(c-1)}}{t^{2/(c-1)}}\right)\right)^c \\
&\leq u+\int_{u}^\infty (u/t)^{2c/(c-1)} dt\\
&= u +\lim_{b\rightarrow \infty}\left(-\frac{c-1}{c+1}u^{2c/(c-1)} t^{-(c+1)/(c-1)} \right)\Big|_u^b \\
&= u+ \frac{c-1}{c+1}u^{2c/(c-1)}u^{-(c+1)/(c-1)}\\
& = u+\frac{1}{2c/(c-1) -1}u = 2c/(c+1)u.
\end{align*}
Therefore, $c_z(e) \leq 2c/(c+1) c_x(e)$ as desired.
\QED

We show that the approximation guarantee for $c$-uniform hypergraphs extend to all hypergraphs $H$ in which each edge contains at most $c$ vertices. In the MRT problem, this corresponds to networks in which each tree edge is covered by at most $c$ switches.

\begin{theorem}\label{thm:cuniform_extension}
Given a network $G=(V,E)$, let $T\subseteq E$ be a spanning tree on $G$ with edge weights $b(e) \in \mathbb{R}_+$ for all $e\in T$. Let the set of switches be $S=E\setminus T$, and let $c= \max_{e\in T} |\{s\in S: T-e+s\text{ is connected}\}|$ be the maximum coverage of any edge in the tree. Then, if $c\geq 2$, there is a polynomial-time $(\frac{2c}{c+1})^2$-approximation algorithm for the MRT problem.
\end{theorem}
\emph{Proof.} For each edge $e$ with fewer than $c$ vertices, modify $H$ by adding additional dummy vertices belonging only to $e$, until $e$ has exactly $c$ vertices. Solve the LP on the resulting $c$-uniform hypergraph $H'$. The objective value of the LP on $H'$ is the same as the objective value of the LP on the original hypergraph $H$, since assigning all weight given to a dummy vertex at a given time to a non-dummy vertex which covers the same edge cannot increase the LP objective. To obtain a solution $\sigma'$, apply $\alpha$-point rounding with kernel \[
K(t,t')= \begin{cases}
4 \frac{t'(t'+1)}{t(t+1)(t+2)} \cdot \mathds{1}[t\geq t'], c=2,\\
\frac{2c}{c+1} \frac{t'^{2/(c-1)}}{\sum_{i=1}^{t} i^{2/(c-1)}} \cdot \mathds{1}[t\geq t'], c\geq 3,
\end{cases}
\]
to the LP solution on $H'$. By Theorem~\ref{thm:alphaallc}, or Theorem 22 of~\cite{Bansal2021} if $c=2$, the expected objective value of $\sigma'$ is at most $(\frac{2c}{c+1})^2$ times the LP objective. 

Finally, convert $\sigma'$ into an ordering $\sigma$ of the vertices of $H$ as follows. For each dummy vertex in $\sigma'$, exchange its place in the ordering with a non-dummy vertex later in the ordering which covers the same edge. When no more such exchanges can be made, remove the dummy vertex from the ordering. This procedure cannot increase the time at which each edge $e$ is covered, and hence cannot increase the objective. Therefore, the expected objective value of the resulting ordering $\sigma$ of the vertex set $V$ is at most  $(\frac{2c}{c+1})^2$ times the LP objective. 
\QED

We further show that this approximation ratio is tight with respect to the natural LP relaxation, thus generalizing Lemmas 26 and 27 of \cite{Bansal2021}. 

\begin{theorem}\label{thm:alphatight}
The integrality gap of the LP for MSSC on $c$-uniform hypergraphs ($c\geq 2$) is at least $(\frac{2c}{c+1})^2$. Therefore, if the LP relaxation is used as a lower bound for any approximation algorithm, one cannot hope to do better than $(\frac{2c}{c+1})^2$. 
\end{theorem}

We include a detailed proof in Appendix~\ref{app:integralitygap}. The key idea of the proof is to construct a $c$-uniform instance out of the union of disjoint complete $c$-uniform hypergraphs of varying sizes. The LP relaxation can efficiently cover each hypergraph by assigning each vertex weight $1/c$, thus covering each complete $c$-uniform hypergraph on $n$ vertices in time $n/c$. In contrast, any integral solution must select all but $(c-1)$ vertices in each of the complete hypergraphs, which takes $n-c+1$ time steps. By constructing instances from increasingly large numbers of increasingly large complete $c$-uniform hypergraphs, the ratio between the objective values of these two solutions approaches $(2c/(c+1))^2$. This ratio matches the approximation guarantee of Theorem~\ref{thm:alphaallc}, and approaches the general approximation ratio of 4 as $c\rightarrow \infty$. These guarantees are worst-case. Computationally (Section \ref{sec:alpha-point-comp}) we find that the greedy algorithm and the $\alpha$-point rounding approaches are competitive, usually within 5\% of the optimal solution as found by solving the IP.

\section{Balancing multiple reliability and energy loss metrics}\label{sec:local}
Different choices of the initial spanning tree $T\subseteq E$ and reconnection order $\sigma$ on $S=E\setminus T$ can have widely diverging performance on {\sc R-Time} and SAIDI, as well as energy costs dependent on the choice of the tree (see, for example, the system in Figure \ref{fig:spokewheel}). Moreover, as we discuss in Section~\ref{sec:greenboro}, the choice of optimization metric will affect the expected reconnection times of various types of buses (e.g.\ in residential or commercial area) differently.

In this section, we give a local search method for choosing a spanning tree $T\subseteq E$. The idea of local search for maximizing monotone submodular functions goes back to~\cite{Calinescu2011}. Our objectives require us to minimize a {\it product} of sums of supermodular functions over changing network topologies, thus making analysis of the process more difficult. 

\subsection{Branch exchange local search}
Since we would like to simultaneously minimize all three objectives as best as possible to decrease loss of energy and improve service reliability across all customers, our goal is to find a single choice of $T$ and $\sigma$ that achieves strong performance on all three objectives. This becomes a multi-objective optimization problem, which we model as a composite objective
\[
\text{SAIDI} \,\cdot\, \text{{\sc R-Time}} \,\cdot\, \text{Energy} = \frac{\sum_{e\in T}f(e)p(e)t(e)}{\sum_{v\in V} w(v)}\cdot \frac{\sum_{e \in T} p(e)t(e)}{\sum_{e \in T} p(e)}\cdot \sum_{e \in T} r(e) f(e)^2,
\]
then apply the well-known spanning tree local search heuristic \textit{branch exchange}.

The branch exchange search 
takes as input a network $G$ with initial spanning tree $T_0$ and an objective function $F(T,\sigma)$, which can be a combination of multiple metrics on $T$ (such as energy) and $\sigma$ (such as {\sc R-Time} and SAIDI)~\cite{civanlar1988distribution}. In each iteration, the local search tests a random pair $(e,s)\in T\times S$ that is feasible (i.e., $T'=T-e+s$ is still connected). The algorithm then generates an order $\sigma'$ on the new set of switches $S+e-s$ and compares the new objective value $F(T',\sigma')$ to the current value $F(T,\sigma)$. If the objective is improved, the exchange is accepted; otherwise, sample another pair and repeat, until there is no improving step. The component of the composite product objective that depends on $\sigma'$, \text{SAIDI $\cdot$ {\sc R-Time}}, can no longer be written an instance of MRT.   Therefore, we use a version of the greedy algorithm to choose $\sigma'$, where at each step we choose the switch $s$ achieving the greatest marginal increase in the product of the weights $p(e)f(e)$ and $p(e)$ used for SAIDI and {\sc R-Time}, respectively.

A single local search step may check dozens or hundreds of possible branch exchanges, and for each possible exchange, the connectivity information needed for computing reconnection orders changes. Over networks of practical interest, existing methods for detecting cycles are very slow. In Appendix~\ref{app:coverage}, we give a method for updating connectivity which achieves a factor of $\Omega(\min\{|E-T|,|V|/d\})$ speedup, where $d$ is the length of the longest cycle in $G$.

If the objective function is chosen to be linear in the weights $b(e)$ of the edges, for instance a linear combination of the {\sc R-Time} and SAIDI objectives, then it is possible to find an optimal spanning tree and reconnection order with an integer programming formulation (Appendix~\ref{sec:app_IP}). However, this program is large and solving it is intractable for larger networks.

\section{Computations}\label{sec:computations}
In this section, we test our methods for improving outage metrics on the Greensboro, NC urban-suburban synthetic network from the NREL SMART-DS project to evaluate the performance of our $\alpha$-point rounding and local search methods. All our code is available via \href{github}{github}
, where we have included the relevant data from NREL SMART-DS dataset in an accessible CSV format.

\subsection{Greensboro dataset}\label{sec:greenboro}

The NREL SMART-DS synthetic network covers most of the Greensboro area, including residential, industrial, and other areas, and contains 145052 buses~\cite{PALMINTIER2021106665}. The network contains demand data for each bus, which we treat as the vertex weights $w(v)$, and lines between buses. We take the failure probability $p(e)$ of each line $e$ to be proportional to the straight-line distance between its endpoints, as in practice, failure rates tend to be proportional to length~\cite{bolacell2020}. The network contains 19 distinct connected components. Of these, 18 correspond to radial distribution networks, each with its own substation, while the other links the substations. We apply our methods to each of the 18 distribution networks, treating them as entirely separate and independent networks.

Some components of the Greensboro network have over 10,000 vertices each. Before performing experiments, we preprocess the network to simplify and reduce its size while retaining key information about its structure. For instance, we remove low-demand degree-1 and degree-2 nodes by contracting incident edges and distributing their demands to neighboring vertices. These steps eliminate leaves and bridge vertices corresponding to buses with low weight. See Appendix~\ref{app:preprocessing}. Repeatedly applying these steps significantly reduces the size of the network. Using a threshold of $W=10~\text{kW}$ on the Greensboro network, the resulting contracted trees have an average of 17\% of the the number of nodes of the original trees. These simplifications primarily consist of contracting leaves corresponding to low-demand buses in residential areas, and do not significantly change the network topology.

\subsection{Adding switches}\label{sec:switches}

\begin{figure}[t]
    \centering
    \includegraphics[width=.49\textwidth]{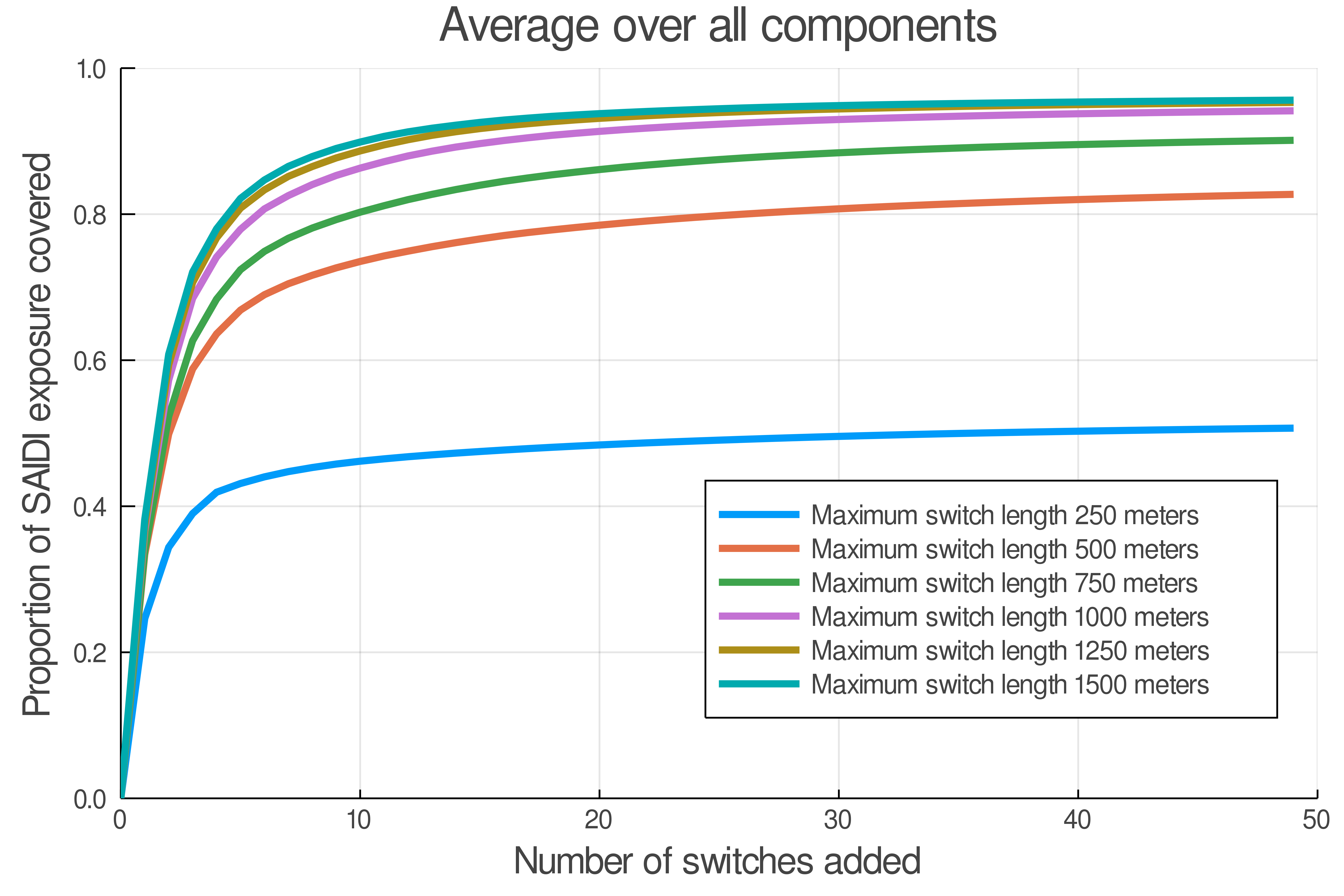}
    \includegraphics[width=.49\textwidth]{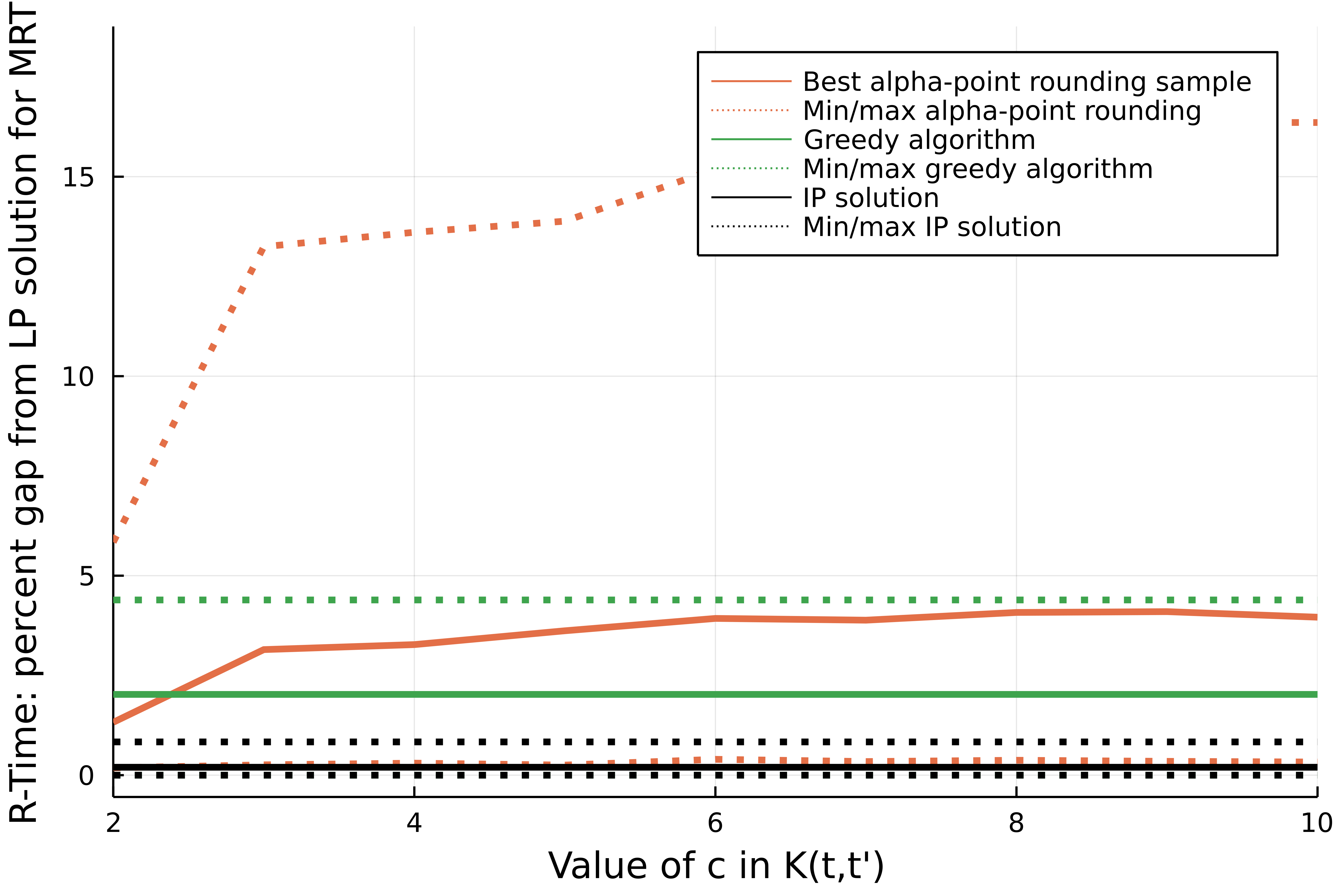}
    \caption{({\bf Left}): the fraction $\sum_{e \text{ covered by } s \in S} x(e)/\sum_{e \in T} x(e)$ of the total exposure covered by the switches added by the greedy algorithm for switch addition, averaged over all components. ({\bf Right}): Performance of six methods for minimizing {\sc R-Time} over the 18 components of the Greensboro network, compared to the LP objective. Solid lines denote average performance over the 18 components, while dotted lines indicate the largest and smallest differences.}
    \label{fig:alpha-point_MRT}
    \label{fig:greedy_coverage}
\end{figure}

The synthetic network has only 108 switches (non-spanning tree edges). To increase coverage, we introduce additional switches, aiming to maximize coverage while minimizing the number and length of the added switches. To do so, we use a two-step procedure for each contracted component. We generate a large set $R$ of candidate switches of length less than some $\ell$. One can solve a redundancy allocation problem to find an optimal set of additional switches $S\subseteq R$, for instance~\cite{park2020milp}, but for our purposes, a simple greedy algorithm suffices to choose switches that provide good coverage while being inexpensive to construct. Define the \textit{exposure} of each edge $e\in T$ to be $x(e) = f(e)p(e)$, its coefficient in SAIDI. At each step, the algorithm adds the switch $s$ that maximizes the weighted coverage $\sum_{e\text{ covered by $s$}} 2^{-|\{t \in S: t \text{ covers } e\}|} x(e)$. For most components, excellent coverage (over 90\% of total exposure) can be achieved by setting $\ell=1000$m and adding at most 20 switches from $S$. Figure~\ref{fig:greedy_coverage} shows how coverage improves as the number of switches added increases, and how this depends on the choice of $\ell$. By adding a number of switches equal to an average of 0.18\% of the total number of tree edges in each component, we are able to cover 90\% of SAIDI exposure. 

\subsection{Approximations for MRT}\label{sec:alpha-point-comp}
Having computed switch sets $S$ for each component network using a maximum length of 1000m, we examine reconnection orders for {\sc R-Time} and SAIDI. To do so, we compare three different approaches for solving MRT: (i) greedy algorithm, (ii) $\alpha$-point rounding on a linear program solution, and (iii) the integer program. Greedy is extremely computationally inexpensive (hence our use of it in our local search algorithms) and solving the LP takes the majority of the time for $\alpha$-point rounding. On many components, solving the IP was of comparable speed to $\alpha$-point rounding, but on some it took much longer. 

Table~\ref{tab:MRTvsSAIDI} gives statistics on the 18 distribution networks and their outage metrics when using the IP for MRT to minimize {\sc R-Time} and SAIDI in turn.
For the reconnection orders $\sigma$ corresponding to each objective, we compute the values of both metrics. The choice of objective matters substantially. The value of SAIDI is on average more than 22\% greater when the MRT objective is {\sc R-Time} than when its objective is SAIDI, and the value of {\sc R-Time} is on average almost 10\% greater when SAIDI is the objective than when {\sc R-Time} is the objective.

Figure~\ref{fig:alpha-point_MRT} compares the performance of the three methods, expressed as a percentage difference relative to the LP objective. For each method, the solid line shows the average performance across the 18 components, and the dotted lines show the best and worst performance over all components. For $\alpha$-point rounding, 500 sets of samples of rounding points $\alpha_s$ were taken for each component and each value of $c$ in $K(t,t') = 2c/(c-1)(t(t+1)^{2/(c-1)})$, and the sample yielding the best objective value was chosen. Note that $\alpha$-point rounding performs best for small values of $c$ where the kernel decays more aggressively, resulting in a solution that tends to be closer to the LP solution, even when some edges in the tree are covered by many switches and the applicable approximation guarantees do not hold. Overall, both polynomial-time methods are competitive, usually within a few percentage points of the optimum, with the greedy method sometimes obtaining the optimal solution and sometimes being outperformed by $\alpha$-point rounding.

\subsection{Choice of optimization metric}

Of the 145052 buses in the network, 63 run at a medium voltage (MV) with a demand of at least 700 kW, while the rest run at a low voltage (LV) with demand of at most 350 kW. The MV buses correspond to commercial and industrial areas, while most LV buses occur in residential areas. Minimize SAIDI will tend to prioritize outages affecting the higher-demand MV buses, whereas minimizing MRT will prioritize restoring power to the more numerous LV buses. Indeed, optimizing for SAIDI produces an 8\% longer expected outage time for LV buses compared to optimizing for {\sc R-Time} (2.70 units, when optimizing for {\sc R-Time} and 2.91 units when optimizing for SAIDI), but a 33\% reduction in expected outage time for MV buses (3.86 units when optimizing for {\sc R-Time} and 2.57 units when optimizing for SAIDI). 

\subsection{Local search for multi-objective optimization}
\begin{figure}[t]
    \centering
    \includegraphics[width=.49\textwidth]{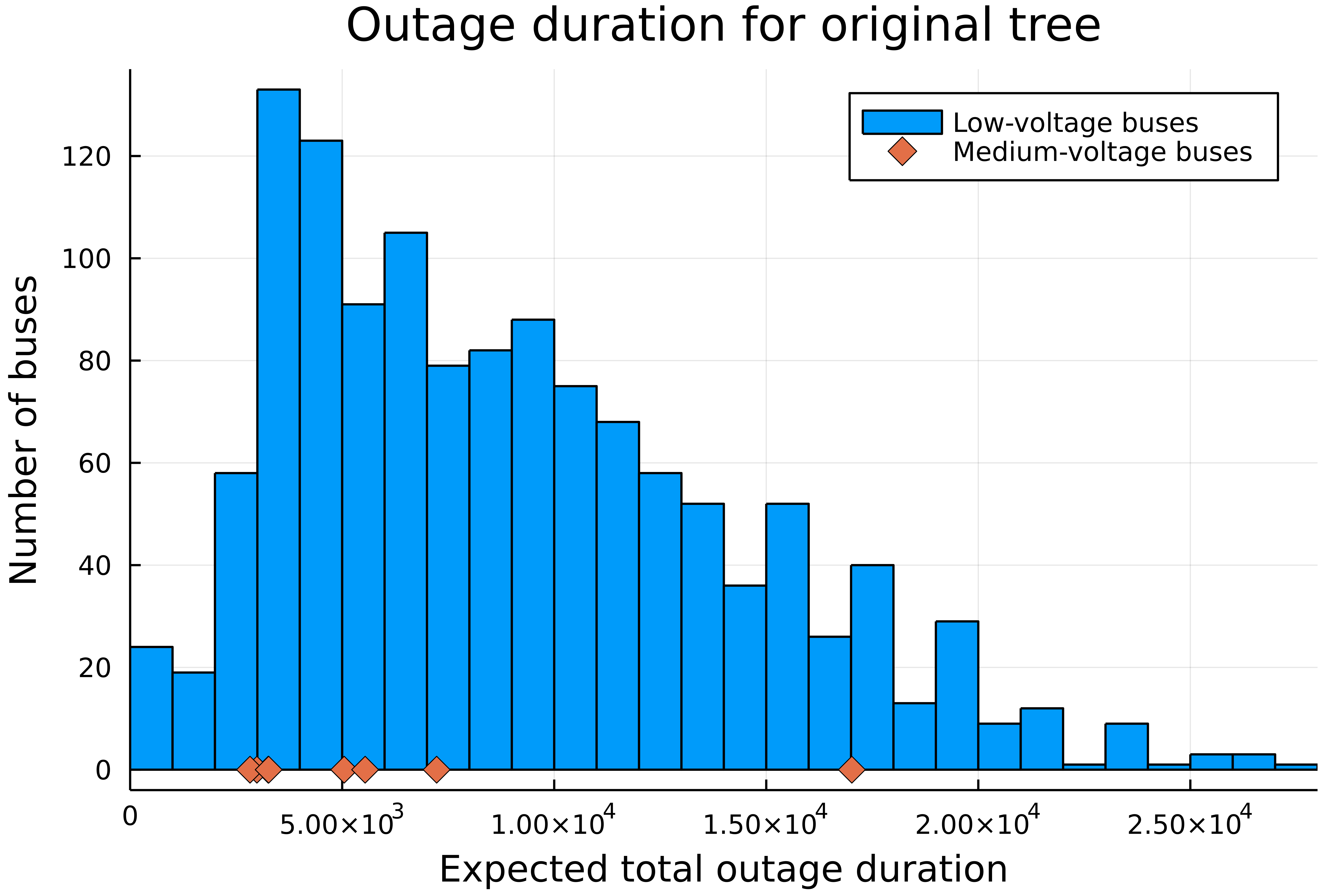}
    \includegraphics[width=.49\textwidth]{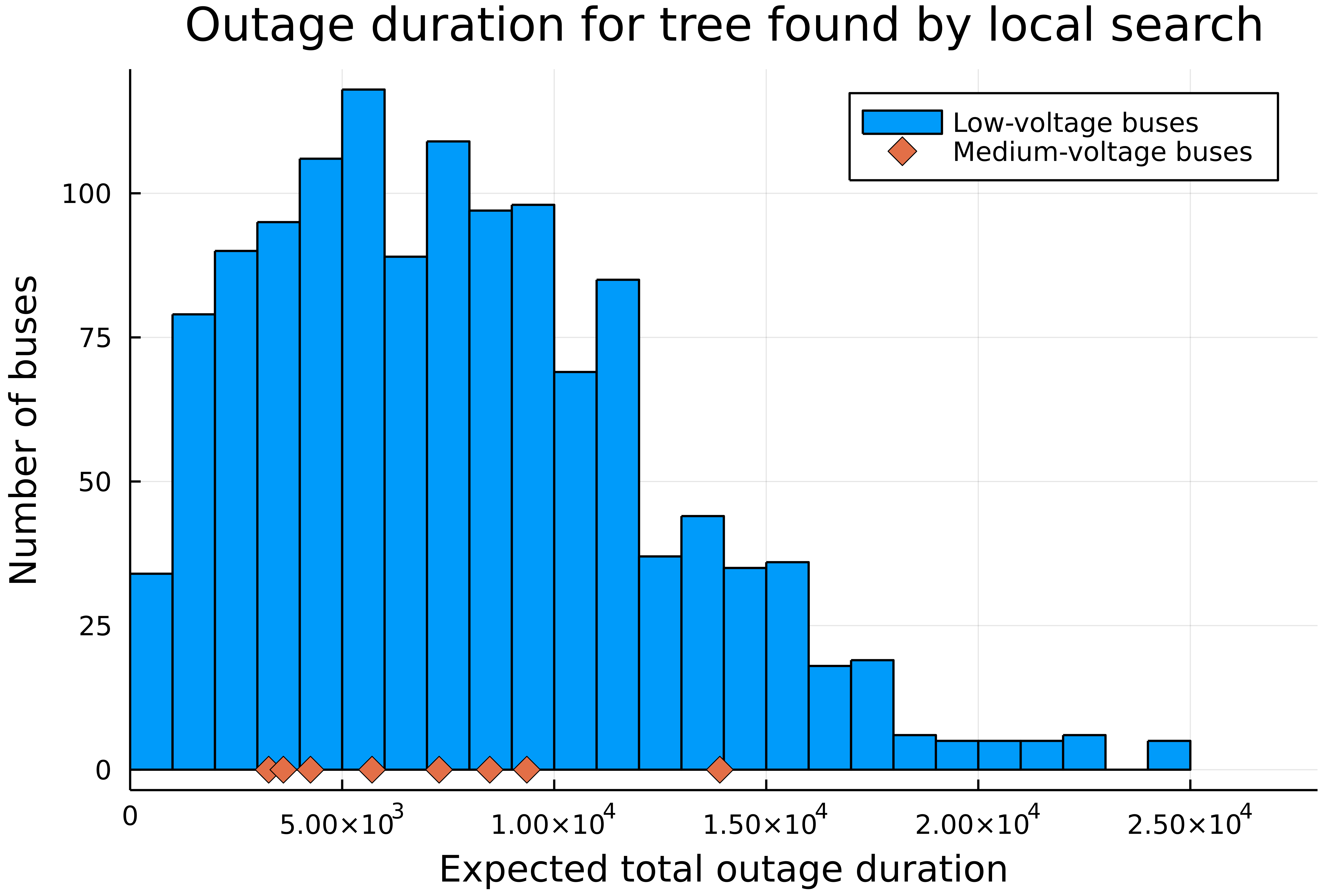}
    \caption{Outage durations for low voltage and medium voltage vertices in component 18. The outage lengths for the 1290 low voltage vertices are plotted as a histogram, while the outage lengths for the 8 medium voltage vertices are plotted individually. (Left): Outage durations for the contracted tree $T$, after adding switches to $G$. (Right): Outage durations for the tree obtained by local search after performing 48 branch exchange steps on $T$. After local search, the distribution of expected outage duration has a higher mass towards lower values for both the groups.}
    \label{fig:localsearch7}
    \vspace{-0.2cm}
\end{figure}

In the branch exchange for local search, for any potential exchange of $e\in T$ and $s\in S$, the energy cost is easy to compute. To compute the value of SAIDI $\cdot$ {\sc R-Time} efficiently we  use a greedy approach, where for each time step $i$ we choose the switch $s$ which achieves the greatest marginal increase in the objective SAIDI$\cdot${\sc R-Time}. Figure~\ref{fig:localsearch7} (left) shows the results of 25 different applications of local search to each component of the Greensboro network, each stopping after 100 exchanges or when a local optimum was reached. For each component, we take the local search run achieving the best product and express the value of each objective at each step as a ratio to its initial values on that component. Finally, we plot the average of these ratios over time for all 18 components. 

\begin{table}[t]
\label{tab:MRTvsSAIDI}
\centering 
\footnotesize 
    \begin{tabular}{|r|r r r|r r|r r|c c|}
    \hline
    \# & Buses & Existing & Added & \sc{SAIDI*} & \sc{R-Time} & \sc{SAIDI} & \sc{R-Time*} & \sc{SAIDI} &  \sc{R-Time} \\
    &       & switches      & switches & &  &  & & gap & gap \\
    \hline
    \hline
1 & 8529 & 4 & 42 & 189.96 & 384.62 & 216.44 & 364.78 & 13.94\% & 5.44\%\\
2 & 8396 & 7 & 47 & 134.25 & 388.57 & 165.38 & 357.85 & 23.18\% & 8.59\%\\
3 & 13010 & 7 & 44 & 158.12 & 295.18 & 169.23 & 277.53 & 7.03\% & 6.36\%\\
4 & 9841 & 15 & 45 & 193.19 & 284.22 & 210.89 & 263.06 & 9.16\% & 8.04\%\\
5 & 9319 & 1 & 46 & 133.85 & 416.48 & 154.00 & 384.51 & 15.06\% & 8.31\%\\
6 & 6238 & 8 & 43 & 169.52 & 325.24 & 184.27 & 296.57 & 8.70\% & 9.67\%\\
7 & 8538 & 3 & 43 & 161.63 & 383.51 & 192.74 & 344.83 & 19.25\% & 11.22\%\\
8 & 8147 & 8 & 43 & 180.38 & 394.75 & 217.34 & 365.96 & 20.49\% & 7.87\%\\
9 & 7286 & 8 & 44 & 244.34 & 289.80 & 450.27 & 235.67 & 84.29\% & 22.97\%\\
10 & 5901 & 5 & 48 & 99.75 & 333.90 & 109.95 & 314.42 & 10.23\% & 6.20\%\\
11 & 11031 & 8 & 43 & 201.56 & 330.38 & 229.29 & 299.16 & 13.76\% & 10.44\%\\
12 & 11536 & 6 & 44 & 147.21 & 482.29 & 167.51 & 455.58 & 13.79\% & 5.86\%\\
13 & 2427 & 3 & 42 & 186.40 & 354.30 & 222.51 & 291.86 & 19.37\% & 21.39\%\\
14 & 6950 & 3 & 36 & 125.78 & 350.83 & 132.27 & 335.37 & 5.16\% & 4.61\%\\
15 & 5914 & 4 & 42 & 116.97 & 371.88 & 121.81 & 345.20 & 4.14\% & 7.73\%\\
16 & 7006 & 2 & 44 & 178.91 & 512.01 & 328.14 & 466.29 & 83.41\% & 9.81\%\\
17 & 4008 & 2 & 44 & 92.37 & 531.91 & 120.52 & 470.12 & 30.47\% & 13.14\%\\
18 & 10953 & 14 & 46 & 240.52 & 367.17 & 303.93 & 330.78 & 26.37\% & 11.00\%\\

\hline
\end{tabular}
\caption{Graphical data and outage metrics for the 18 distribution networks. The gaps are the relative increase in objective (the starred metric) from using the integer programming formulation to optimize for the other metric.}
\end{table}

Since the possible branch exchanges are sampled randomly, the results of the local search can vary across multiple runs. Figure~\ref{fig:localsearch7} (right) shows the results of 50 different applications of local search to component 9 of the Greensboro network, each stopping after 90 exchanges. For each step, the median of each metric across all 50 runs is plotted in the band between the 10th and 90th percentile values. All three objectives improved substantially during the local search, with energy losses decreasing by 85\%, SAIDI by 43\%, {\sc R-Time} by 17\%, and their product by 93\%. Similar plots for other components are included in Appendix~\ref{sec:app_component_plots}, with some variation in the degree of improvement using local search. 

\begin{figure}[t]
    \centering
    \includegraphics[width=.48\textwidth]{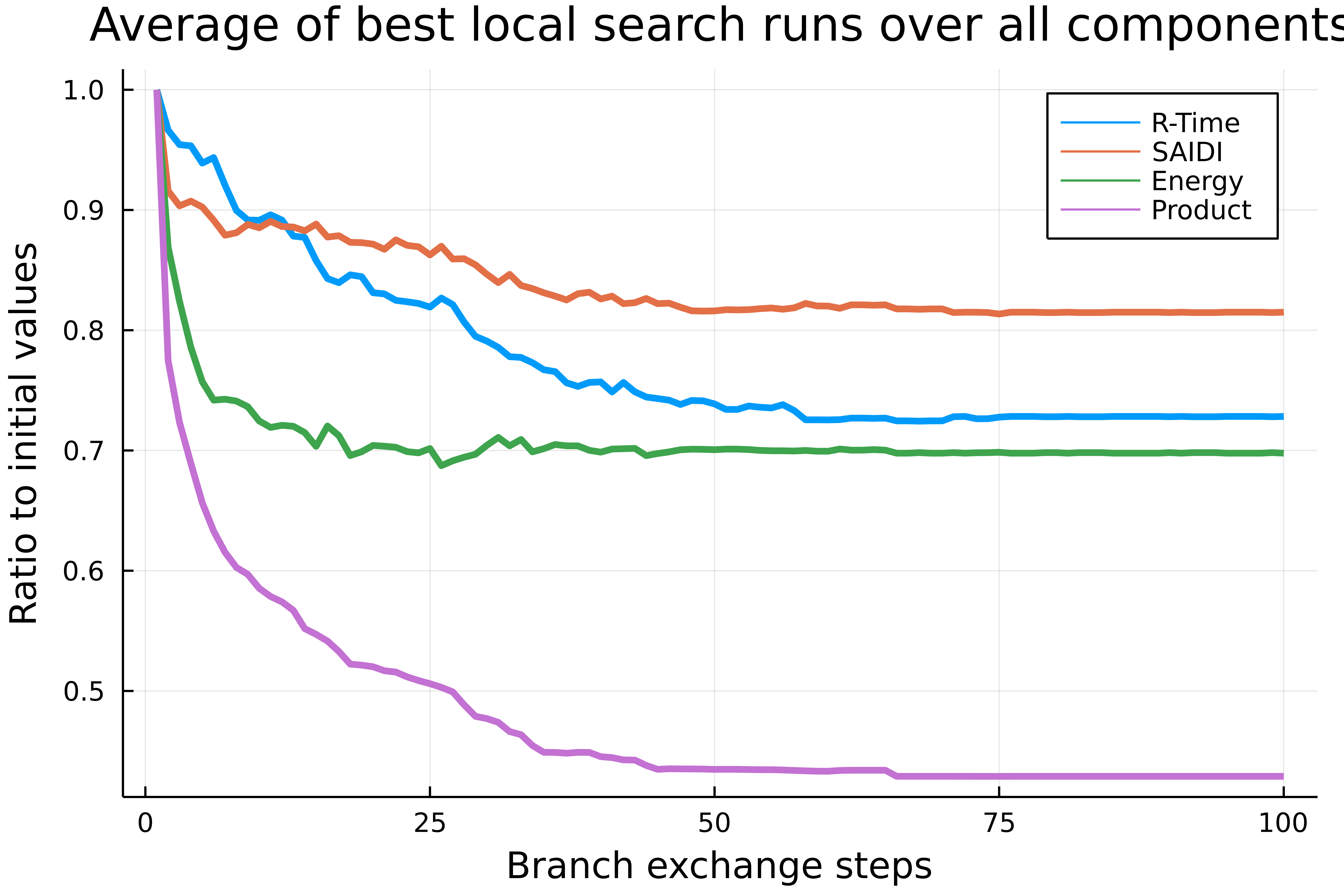}
    \includegraphics[width=.49\textwidth]{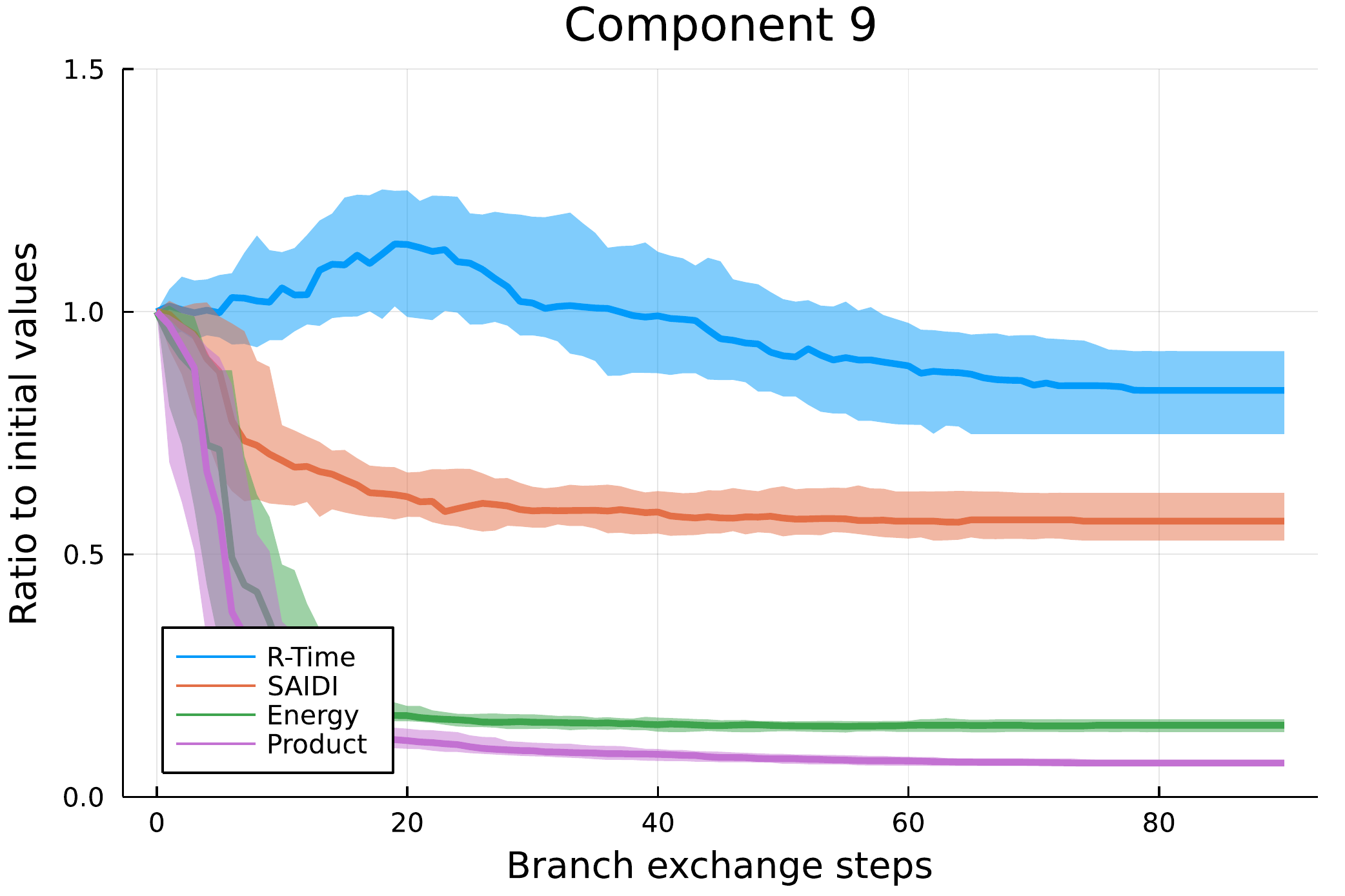}
    \caption{{\sc R-Time}, SAIDI, energy, and product objectives during local search. Objectives are plotted as a ratio to the initial values. \textbf{Left:} Average best values for each objective, taken over 25 applications of local search on all components. \textbf{Right:} Median values and percentiles (first to last deciles), taken over 50 applications of  local search on component 9, highlighting the large improvements made by local search and how {\sc R-Time} and {\sc SAIDI} often behave differently.}
    \label{fig:local_search_equity}
    \vspace{-0.2cm}
\end{figure}

Moreover, we find that as an artifact of the local search, reliability across groups of nodes improves, and the gaps in the expected outage time decrease significantly. Recall that the outage time of a vertex $v$ is 
$\sum_{e:v\in D(e)}p(e)\min_{s\in E\setminus T:\text{ $s$ covers $e$}} \sigma(s)$, where $D(e)$ is the set of edges in the unique path in $T$ between $e$ and the vertex $r$ corresponding to the substation, $p(e)$ is the probability of a fault occurring on edge $e$, and $\sigma$ is the reconnection order selected by the greedy algorithm to minimize the product of {\sc R-Time} and SAIDI. Note that our objectives are not defined with respect to any grouping of nodes by demographic or voltage size\footnote{This may be of interest to the fair OR community, since we are able to move away from the ``group-fairness'' notions, and do not explicitly use groups in the objective.}. However, we find that the difference in expected outage time for medium and low voltage nodes (i.e., industrial and residential nodes, respectively) decreases significantly. For instance, in component 18, which has the most medium voltage nodes (eight) of any component, low voltage nodes had an expected outage time on average 53.3\% greater than that of medium voltage nodes, and modifying the base tree using local search reduced this gap to 12.3\%. This improvement in service equity is evident from the histogram of expected outage times for all the nodes in the network (see Figure \ref{fig:local_search_equity}). 

\section{Conclusion}
We analyzed a decentralized self-healing approach for automatically reconfiguring a distribution system into an operational radial network after a fault occurs by creating an ordering in which switches automatically close upon detection of a downstream fault. The choice of the switches' waiting times, which requires small intervals between each switch's assigned times to avoid creating cycles in the network, is equivalent to choosing a reconnection ordering. The choice of this ordering significantly impacts the expected time to reconnect under normal disruptions and thus affects reliability metrics such as SAIDI, which are the basis for regulator-imposed financial incentives for performance.

We modeled SAIDI and expected reconnection time {\sc R-Time} using a new problem called Min Reconnection Time (MRT), which we showed is a special case of the well-known minimum linear ordering problem. We showed that MRT is NP-hard and generalized the kernel-based rounding approaches of Bansal et al.\ to interpolate approximation guarantees from 9/4 to 4 dependent on the coverage of each tree edge. Finally, using local search, we optimized multiple metrics simultaneously on the NREL SMART-DS Greensboro dataset, thereby giving a proof of concept that branch exchange might be implemented efficiently using proxies for MRT subproblems. 

We demonstrated that the choice of outage metric during the optimization process has a large impact on the behavior of the reconnection order. Optimizing for SAIDI resulted in a 33\% reduction in reconnection time for the medium voltage buses that occur in commercial and industrial areas, while optimizing for
{\sc R-Time} reduced reconnection time in the more numerous but lower-demand low voltage buses. Depending on the priorities of the network designers, incorporating both metrics into the optimization process, as we do in our branch exchange local search, may be appropriate.

Our work raises interesting theoretical and practical open questions. First, it is open if the greedy algorithm for MRT attains a worst-case ratio better than 4-factor, as our experiments suggest. Second, local search is currently well-understood for maximization of submodular functions under matroid base constraints~\cite{Calinescu2011}, but our  local search minimizes the {\it product} of supermodular functions over the matroid base constraint. This seems to be a non-trivial extension, and we leave it as an open question. Lastly, we hope that our work provides a starting point for analyzing multiple simultaneous failures in the network, a situation that might be interesting for major disruptions in practice.

\bibliographystyle{spmpsci} 
\bibliography{main.bib}

\clearpage
\begin{appendix}

\section{Motivation from electrical distribution networks}\label{app:elec_networks}

Our graph-based model is a simplified representation of the distribution network compared to the model used in, for instance, power flow studies that additionally involve information about the line reactances, reactive demands, capacitor susceptances, etc.~\cite{kersting2017}. While precluding details such as reactive power injections and voltages, our approach enables application of the graph theoretic concepts that are the basis for this paper's methodology.\footnote{Similar models have been used in power systems work such as \cite{kersulis2018pscc,Khodabakhsh2017,li2020}.}

In the network $G=(V,E)$, vertices $v\in V$ correspond to buses and edges $e\in E$ correspond to lines. The initial configuration of the network is a spanning tree $T$, with the edges that are initially inactive corresponding to the set of switches $S = E\setminus T$. We specify a reconnection order $\sigma: S\longrightarrow [|S|]$ on the switches. Events which cause tree edges to become inactive correspond to electrical faults. Reconnecting the network by activating an edge corresponds to closing a switch on a line not in $T$, thus energizing it. The constraint of energized lines forming a radial structure is common in distribution networks~\cite{lavorato2012,wang2020radial,kersting2017}. It is particularly crucial to maintain this structure after a fault, as otherwise we risk energizing the fault.

Some edges have breakers which will trip and stop flow in the event of a nearby fault, thus isolating it even after the network is reconnected and power is restored to vertices that had lost it. In general, not every edge will have a breaker, but if necessary we can preprocess the graph by contracting each set of connected edges in $T$ that has only one breaker to a single edge $e'$. We set the probability $p(e')$ that $e'$ experiences a fault in a given time period to be the sum of the fault probabilities $p(e)$ of the contracted edges. 

After a fault, every bus receiving flow will have positive voltage, and the others will have zero voltage; the sensors and relays employed on switches $s$ read the voltages of the incident buses to determine whether they are receiving power and thus decide whether to wait $\sigma(s)$ time before closing if one bus still has zero voltage. Such relays may be expensive, but installing them on a small number of lines is feasible for current distribution systems~\cite{byeon2020communication}.

\section{Proof of Lemma~\ref{lem:isum}}\label{app:isumlemmaproof}

We use Lemma~\ref{lem:isum} as a step in the proof of Theorem~\ref{thm:alphaallc}, and this lemma in turn requires Lemma~\ref{lem:HHlike}, as stated in this section.  

{\sc Lemma~\ref{lem:isum}.} For $t\in\mathbb{N}$ and $0<p\leq 1$, $\sum_{i=1}^{t} i^p \leq\frac{p}{p+1}\frac{t^p (t+1)^p}{(t+1)^p-t^p}$.

\emph{Proof.} We proceed by induction on $t$. The claim holds for $t=1$ since $2^p-1\leq p$ and therefore
\[
\sum_{i=1}^{t} i^p = 1 \leq \frac{p\cdot 2^p}{(p+1)(2^p-1)} = \frac{p}{p+1}\frac{t^p (t+1)^p}{(t+1)^p-t^p}.
\]
Now suppose $\sum_{i=1}^{t} i^p \leq\frac{p}{p+1} \frac{t^p (t+1)^p}{(t+1)^p-t^p}$ holds for some $t\geq 1$. Then
\begin{align*}
\sum_{i=1}^{t+1} i^p &\leq \frac{p}{p+1}\frac{t^p (t+1)^p}{(t+1)^p-t^p} + (t+1)^p\\
&= \frac{p}{p+1}\frac{t^p (t+1)^p + \frac{p+1}{p}(t+1)^{2p}-\frac{p+1}{p}(t+1)^p t^p}{(t+1)^p-t^p} \\
&= \frac{p}{p+1}(t+1)^p\frac{\frac{p+1}{p}(t+1)^{p}-\frac{1}{p}t^p}{(t+1)^p-t^p}\\
&= \frac{p}{p+1}(t+1)^p\left(\frac{(t+1)^{p}}{(t+1)^p-t^p}+\frac{1}{p}\right). \numberthis\label{eqn:isum}
\end{align*}
The derivative of $f(t)=\frac{(t+1)^{p}}{(t+1)^p-t^p}$ is
\[
f'(t) = p\frac{(t+1)^{p}t^{p-1}-(t+1)^{p-1}t^{p}}{((t+1)^p-t^p)^2} = p\frac{(t+1)^{p}t^{p-1}-(t+1)^{p-1}t^{p}}{\left(p\int_t^{t+1} x^{p-1} dx\right)^2}. 
\]
To bound the denominator of the derivative, we use the following lemma, due to~\cite{Pinelis-MO}.\footnote{This lemma is similar to the bound of $(v-u) (u^{q}+v^{q})/2$ given by the Hermite-Hadamard inequality on convex functions, but with the arithmetic mean replaced by the geometric mean.}
\begin{lemma}\label{lem:HHlike}
Let $q\in(-1,0]$, $v\geq 1$, and $u\in (0,v)$. Then 
\[
\int_u^v x^q dx \leq (v-u)u^{q/2}v^{q/2}.
\]
\end{lemma}
\emph{Proof.}
The desired inequality is homogenous of degree $q+1$ in $u$ and $v$. Therefore, by rescaling $u$ and $v$ simultaneously, we may assume without loss of generality that $u=1$. Then the desired inequality is equivalent to 
\[
(q+1)\left(\int_u^v x^q dx - (v-u)u^{q/2}v^{q/2}\right) = v^{q+1}-1-(q+1)(v-1)v^{q/2} \leq 0.
\]
We next show that for all $q$, the derivative with respect to $v$ of the expression $g(v,q)=v^{q+1}-1-(q+1)(v-1)v^{q/2}$, 
\begin{align*}
\frac{d}{dv}g(v,q)&=(q+1)v^q - (q+1)\left(v^{q/2}+\frac{q}{2}(v-1)v^{q/2-1}\right), \\
&=\frac{q+1}{2}v^{q/2-1}(2v^{q/2+1}-2v-q(v-1)),
\end{align*}
is negative.

As a function of $q$ for any fixed $v$, $h_v(q)=2v^{q/2+1}-2v-q(v-1)$ is convex. Moreover,  $h_v(-1)=-(\sqrt{v}-1)^2<0$ and $h_v(0)=0$. Therefore $h_v$ is negative for all $q\in (-1,0]$, and since $\frac{q+1}{2}v^{q/2-1}\geq0$, the derivative $\frac{d}{dv}g(v,q)$ is negative as desired. For all $q$ and $v=1$, $g(v,q)=0$, so since $g(v,q)$ is decreasing in $v$, $g(v,q)\leq 0$ for all $q$ and $v$. This proves the desired inequality.
\QED

For all $u>0$, applying Lemma~\ref{lem:HHlike} with $q=p-1$ and $v=u+1$ yields 
\begin{align*}
f'(u)&=p\frac{(u+1)^{p}u^{p-1}-(u+1)^{p-1}u^{p}}{\left(p\int_u^{u+1} x^{p-1} dx\right)^2}, \\
&\geq p\frac{(u+1)^{p}u^{p-1}-(u+1)^{p-1}u^{p}}{p^2(u+1)^{p-1}u^{p-1}},\\
&= \frac{(u+1)-u}{p} = \frac{1}{p}.
\end{align*}
Therefore $f(t+1) \geq f(t)+1/p$, and from Equation~\ref{eqn:isum} we have that
\begin{align*}
\sum_{i=1}^{t+1} i^p &\leq \frac{p}{p+1}(t+1)^p\left(\frac{(t+1)^{p}}{(t+1)^p-t^p}+\frac{1}{p}\right), \\
&= \frac{p}{p+1}(t+1)^p (f(t)+1/p), \\
&\leq  \frac{p}{p+1}(t+1)^p f(t+1),  \\
&= \frac{p}{p+1}\frac{(t+1)^p (t+2)^p}{(t+2)^p-(t+1)^p}, 
\end{align*}
which completes the induction.
\QED

\section{Proof of the Integrality Gap in  Theorem~\ref{thm:alphatight}}\label{app:integralitygap}
\emph{Proof.}
For a suitable $k=\omega(1)$ and $i=1,\dots,k$, let $n_i=Ni^{-\alpha}$ where $\alpha = 2/(c+1)+\varepsilon$, with $\varepsilon$ approaching 0 and $N$ large enough that $n_i$ can be rounded without affecting the solution. Let $H$ be the hypergraph consisting of disjoint copies of the complete $c$-uniform hypergraphs $K_i$ on $n_i$ vertices.

We first upper-bound the LP cost of MSSC on $H$. Setting
\[
x_{v,t} = \begin{cases} 
\frac{1}{n_i},\text{ if } v\in K_i \text{ and } \sum_{j<i}\frac{n_j}{c} <t \leq \sum_{j\leq i}\frac{n_j}{c}\\
0,\text{ otherwise,}\end{cases}
\]
gives a feasible solution. In this solution, each edge $e$ in $K_i$ is completely covered by time $(n_1+\dots+n_i)/c$, and so the LP objective is at most
\begin{align*}
    \sum_{i=1}^{k}\binom{n_i}{c}\frac{n_1+\dots+n_i}{c} & \leq \sum_{i=1}^{k}\frac{N^{c}}{c!} i^{-c\alpha} \frac{N}{c} \sum_{j=1}^{i} j^{-\alpha}\\
    &\leq \sum_{i=1}^{k}\frac{N^{c+1}}{c!c} i^{-c\alpha}  \frac{i^{1-\alpha}}{1-\alpha}\\
    & \leq\frac{N^{c+1}}{c!c(1-\alpha)}  \sum_{i=1}^{k} i^{1-(c+1)\alpha} = \frac{N^{c+1}}{c!c(1-\alpha)}  \sum_{i=1}^{k} i^{-1-(c+1)\varepsilon}. 
\end{align*}

As we make $\varepsilon$ arbitrarily small, $\sum_{i=1}^{k} i^{-1-(c+1)\varepsilon} \leq 1+\int_{1}^{\infty}x^{1-e\varepsilon}\ dx = 1+1/(c+1)\varepsilon.$ Thus, the objective of the LP solution is at most 
\[
\frac{N^{c+1}}{c!c(1-\alpha)}  \sum_{i=1}^{k} i^{-1-(c+1)\varepsilon} \leq \frac{N^{c+1}}{c!c(\frac{c-1}{c+1}-\varepsilon)} \frac{1}{(c+1)\varepsilon}.
\]

This gives a bound of approximately $N^{c+1}/c!c(c-1)\varepsilon$ for the LP solution for sufficiently small $\varepsilon$ and large $k$.

Next, we show that the greedy algorithm gives an optimum solution to MSSC on $H$, and lower-bound the cost of this solution. In any solution $\sigma$ to MSSC, the $j$th vertex selected from $K_i$ covers $\binom{n_i-j+1}{c}$ hyperedges that had not been covered before. Therefore, the order of vertices chosen within each $K_i$ does not matter, and permuting the order in which vertices from each $K_i$ are chosen relative to each other does not change the number of hyperedges those vertices cover. Thus the greedy algorithm, which at each time step picks a vertex from one of the $K_i$ with the greatest number of unpicked vertices, produces an optimal solution.

For any $1\leq i \leq k-1$, consider the $i(n_i-n_{i+1})$ time steps during which the greedy algorithm selects vertices from $K_1,\dots,K_i$ until $n_{i+1}$ vertices remain in each. During each of these time steps, there are $\binom{n_{i+1}}{c}$ uncovered edges in each of $K_1,\dots,K_i$, and $\binom{n_{j}}{c}$ uncovered edges in each $K_j$ for $j\geq i+1$. Therefore, the objective cost is at least
\[
\sum_{i=1}^{k-1} i(n_i-n_{i+1}) \left(i\cdot\binom{n_{i+1}}{c} +\sum_{j=i+1}^{k} \binom{n_{j}}{c}\right).
\]
Since $\binom{n_i}{c}\approx n_i^c/c!$, we have that
\[
i\cdot\binom{n_{i+1}}{c}+\sum_{j=i+1}^{k}\binom{n_{j}}{c} \approx \frac{N^c}{c!} \left(i^{1-c\alpha} + \int_{i}^{k} j^{-c\alpha}\ dj\right) 
= \frac{N^c}{c!} \left(i^{1-c\alpha} + \frac{1}{c\alpha-1}(i^{1-c\alpha} - k^{1-c\alpha})\right).
\]
Using this and $n_i-n_{i+1}\approx N\alpha i^{-1-\alpha}$, the objective is at least
\[
\frac{N^{c+1}}{c!}\sum_{i=1}^{k}i\cdot \alpha i^{-1-\alpha} \left(\frac{c\alpha}{c\alpha-1}i^{1-c\alpha} + \frac{1}{c\alpha-1}\right) 
= \frac{N^{c+1}}{c!}\left(\sum_{i=1}^{k} \frac{c\alpha^2}{c\alpha-1}i^{1-(c+1)\alpha} + \frac{\alpha k^{1-c\alpha}\sum_{i=1}^{k}i^{-\alpha}}{c\alpha-1}\right).
\]
In this expression, as $k$ increases and $\varepsilon$ decreases, the first term goes to 
\[
\frac{N^{c+1}}{(c-1)!}\frac{\alpha^2}{c\alpha-1}\frac{1}{(c+1)\alpha-2} \approx \frac{N^{c+1}}{(c-1)!}\frac{2^2}{(c+1)^2\frac{c-1}{c+1}}\frac{1}{(c+1)\varepsilon} = \frac{4 N^{c+1}}{(c-1)!(c-1)\varepsilon(c+1)^2}.
\]

Combining the two objective bounds yields the desired integrality gap of at least
\[
\frac{\frac{4 N^{c+1}}{(c-1)!(c-1)\varepsilon(c+1)^2}}{\frac{N^{c+1}}{c!c(c-1)\varepsilon}}=\left(\frac{2c}{c+1}\right)^2. 
\]
\QED

\section{Local search and efficient updates}\label{app:coverage}
\begin{algorithm}[!t]\footnotesize
 \caption{\textsc{Branch exchange local search}}\label{alg:local}
 \begin{algorithmic}[1]
\State \textbf{Input:} Network $G=(V,E)$; initial spanning tree $T_0$ over $V$ and order $\sigma_0$ over $S=E\setminus T$;  probabilities of failure $p(e)$ for each edge $e\in T_0$; vertex weights $w(v)$ for each $v\in V$; objective function $F(T,\sigma)$.
\State \textbf{Output:} A spanning tree $T$ over $V$; a reconnection order $\sigma$ over $E-T$.
\State $T\leftarrow T_0$, $\sigma \leftarrow \sigma_0$, $P\leftarrow T\times (E-T)$
\While{$|P|>1$}
\State $(e,s)\sim P$
\State $P\leftarrow P\setminus \{(e,s)\}$
\State $T'\leftarrow T-e+s$
\State Find a reconnection order $\sigma'$ using greedy or other algorithms for metrics used in $F$.
\If{$F(T',\sigma')<F(T,\sigma)$}
\State $T\leftarrow T', P\leftarrow T\times (E-T)$
\EndIf
\EndWhile
\State \Return $T,\sigma$
\end{algorithmic}
\label{alg:mrt_branch}
\end{algorithm}

Given a network $G=(V,E)$, a substantial portion of the running time of the branch exchange local search Algorithm~\ref{alg:local} is updating the connectivity information, i.e.\ for which pairs of switches $s$ and tree edges $e$ the tree $T-e+s$ remains connected. During the local search procedure, $T$ and $S=E\setminus T$ will evolve, requiring connectivity updates at each step. Indeed, updating this information, which is input to the greedy algorithms for MRT, is significantly more computationally expensive than executing those algorithms, computing energy costs, and comparing the objective values to evaluate whether to accept a branch exchange. This information is useful in filtering and efficiently choosing (at random) feasible branch exchange pairs $(e,s)$.

In this section, we give a description of those pairs of tree edges and switches whose connectivity data might change after a branch exchange update. Updating only these pairs improves the computational efficiency of this step significantly: by a factor of $\Omega(\min\{|E-T|,|V|/d\})$, where $d$ is the circumference, or length of a longest cycle in $G$.

Consider a branch exchange update on $T\subseteq E$. This update replaces $T$ with $T' = T- e + s_e$, where $e\in T$, $s_e\in E- T$, and $T'$ is connected. That is, it exchanges the tree edge $e$ for the switch $s_e$, yielding the new spanning tree $T'$. To compute a reconnection order on the new set of switches $E- T'$, whether by the greedy algorithm or an integer program, it is necessary to know the \textit{covering data} of $G$ and $T'$, i.e.\ for which pairs $(f, s)$ of edges $f\in T'$, $s\in E- T'$ the graph $T'- f + s$ is connected. If we have already computed the covering data for $G$ and $T$, to update this list for $T'$ it suffices to consider only those pairs $(f,s)$ excluded by the following lemma.     

\begin{lemma}\label{lem:exchange_equivalence}
Let $G=(V,E)$ be a graph and let $T\subseteq E$ be a spanning tree of $G$. Let $e \in T$, $s_e \in E- T$ be an edge and a switch such that $T^\prime = T - e + s_e$. Then, all covering data for $T$ also hold for the new spanning tree $T^\prime$, except the following pairs $(f,s)$: 
\begin{enumerate}[leftmargin=18pt]
\item[(i)] $f$ is in the unique cycle $C_e$ in $T + s_e$ and the unique cycle $C_s$ in $T+s$ contains $e$.  
\item[(ii)] $f$ is in the unique cycle $C'_e$ in $T'+ e$ and the unique cycle $C'_s$ in $T'+s$ contains $s_e$.
\end{enumerate} 
In other words, $T- f + s$ is connected if and only if $T'- f + s$ is connected, for all other such pairs of edges $(f,s)$.
\end{lemma}

To prove this result, we show a second lemma:

\begin{lemma}\label{lem:exchange_covering}
Let $G=(V,E)$ be a graph and let $T\subseteq E$ be a spanning tree of $G$. Let $s_e,s_f\in S = E-T$ and let $C_e, C_f$ be the unique cycles in $T+ s_e, T+ s_f$, respectively. Let $e \in C_e \cap T$ and $f\in C_f \cap T$. Then if $e \not \in C_f$ or $f \not \in C_e$, $T- e-f + s_e+s_f$ is connected.
\end{lemma}
\emph{Proof.}
By symmetry, it suffices to show the result when $e \not \in C_f$.

By hypothesis, $T- e + s_e$ is a spanning tree on $G$. Let $u$ and $v$ be the two vertices incident to $f$. Then $T-e-f + s_e$ has two connected components, one containing $u$ and one containing $v$. Therefore, it suffices to show that there is a path between $u$ and $v$ in $T- e-f + s_e+s_f$. Indeed, $C_f - f$ is such a path, since $e \not \in C_f$. Thus $T- e_f +s_e+s_f$ is connected.
\QED

Proof of Lemma~\ref{lem:exchange_equivalence}
Suppose neither of the conditions holds. 

If $T - f + s$ is connected, $f\in C_s$. Since condition (i) does not hold, either $e \not \in  C_s$ or $f \not \in C_e$. Therefore, by Lemma~\ref{lem:exchange_covering}, $T-e-f + s_e+s = T' - f+ s$ is connected.

If  $T'- f + s$ is connected, $f \in C'_s$. Since condition (ii) does not hold, either $f \not \in C'_e$ or $s_e \not \in C'_s$. Therefore, by Lemma~\ref{lem:exchange_covering}, $T'- s_e-f + e+s = T- f + s$ is connected.

Therefore $T- f + s$ is connected if and only if $T'- f + s$ is connected.

In general, computing the connectivity of $T- f + s$ for all pairs $(f, s)$ of edges $f\in T$, $s\in E- T$ takes $O(|V|^2|E- T|)$ time, for instance by performing a depth-first search on each of the $|V||E- T|$ graphs $T- f + s$. Let $d$ be the circumference (length of longest cycle) of $G$. Then by updating connectivity only for those pairs excluded in Lemma~\ref{lem:exchange_equivalence}, the time needed to update connectivity for $T' = T- e + s_e$ is reduced to
\begin{align*}
    O(|V|(|E- T|+|V|+|C_e||E- T|+|C'_e||E- T|) &= O(|V|(|V|+d|E- T|)) \\
    &= O(|V|^2 + d|V||E- T|).
\end{align*}

Therefore, updating only pairs $(f,s)$ as described in Lemma~\ref{lem:exchange_equivalence}, we can save an overall factor of $\Omega(\min\{|E-T|,|V|/d\})$ runtime in each step of the local search procedure.

\section{Network preprocessing}\label{app:preprocessing}
We employ two reduction steps in our preprocessing, which reduces the size of the components and improves tractability: (i) we contract degree 2 vertices while splitting the demand evenly across its two neighbors, and (ii) we contract degree 1 vertices if their demands are less than a threshold. These steps eliminate bridge vertices, as well as leaves corresponding to buses with low weight. We give precise pseudocode in Algorithm~\ref{alg:contraction}. Repeatedly applying these steps significantly reduces the size of the network. Using a threshold of $W=10~\text{kW}$ on the Greensboro network, the resulting contracted trees have an average of 17\% of the the number of nodes of the original trees. These simplifications primarily consist of contracting leaves corresponding to low-demand buses in residential areas, and do not significantly change the network topology. 

\begin{algorithm}[t]\footnotesize
 \caption{(\textsc{Tree Contraction}) Reducing tree size}
 \begin{algorithmic}[1]
\State \textbf{Input:} Spanning tree $T$ over a vertex set $V$;  probabilities of failure $p(e)$ for each edge $e\in T$; weights $w(v)$ for each $v\in V$, weight threshold $W$.
\State \textbf{Output:} A tree $T'$.
\State \textbf{Initialize:} $T'=T$.
\While{a vertex was deleted in the previous loop}
\For{$v \in V$}
    \If{$\deg_{T'}(v)=1$ and $w(v)<W$}
    \State Add $w(v)$ to $w(u)$, where $u$ is $v$'s neighbor.
    \State Add $p(uv)\cdot w(v)/w(u)$ to $p(uw)$, where $w$ is $u$'s parent. Delete $v$.
	\EndIf
	\If{$\deg_{T'}(v)=2$}
	\State Add $w(v)/2$ to $w(u)$ and $w(v)$, where $u$ and $w$ are $v$'s neighbors.
    \State Add edge $uw$ with failure probability $p(uw) = p(uv)+p(vw)$. Delete $v$.
	\EndIf
\EndFor
\EndWhile
\State \Return $T'$.
\end{algorithmic}
\label{alg:contraction}
\end{algorithm}

\clearpage
\section{Plots for local search on other Greensboro components}\label{sec:app_component_plots}

\begin{figure}[H]
    \centering
\includegraphics[width=.49\textwidth]{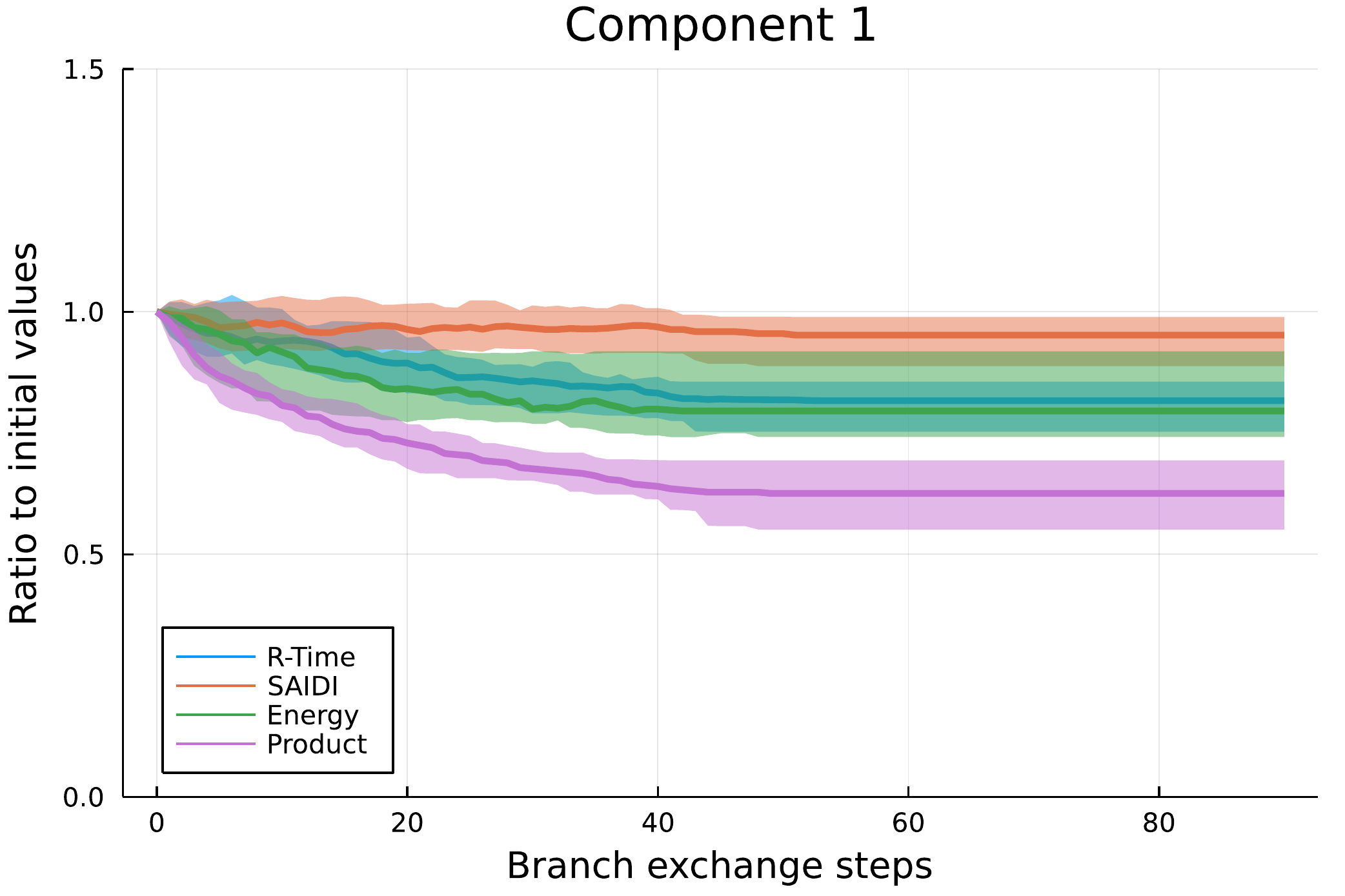}
\includegraphics[width=.49\textwidth]{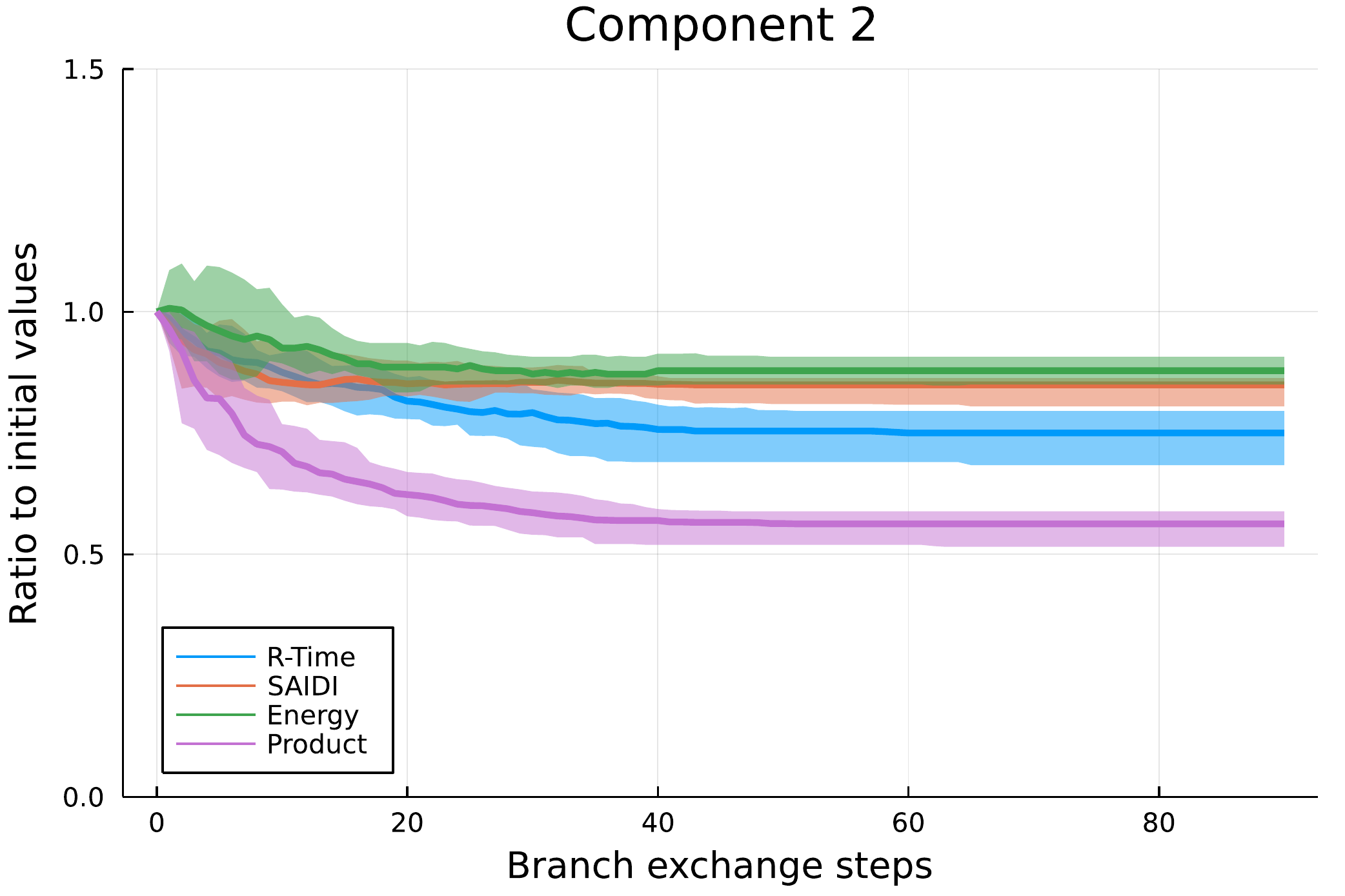}
\includegraphics[width=.49\textwidth]{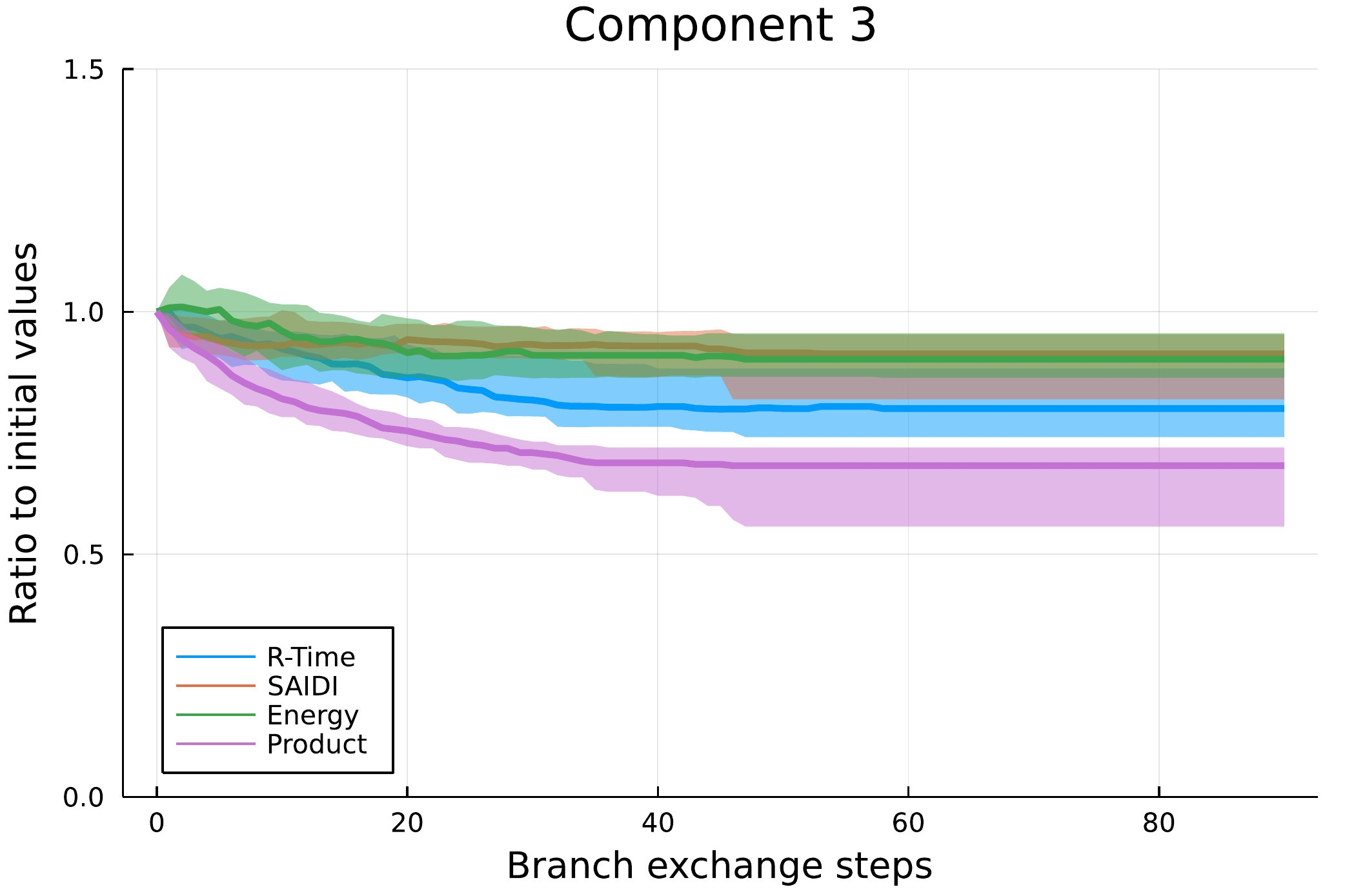}
\includegraphics[width=.49\textwidth]{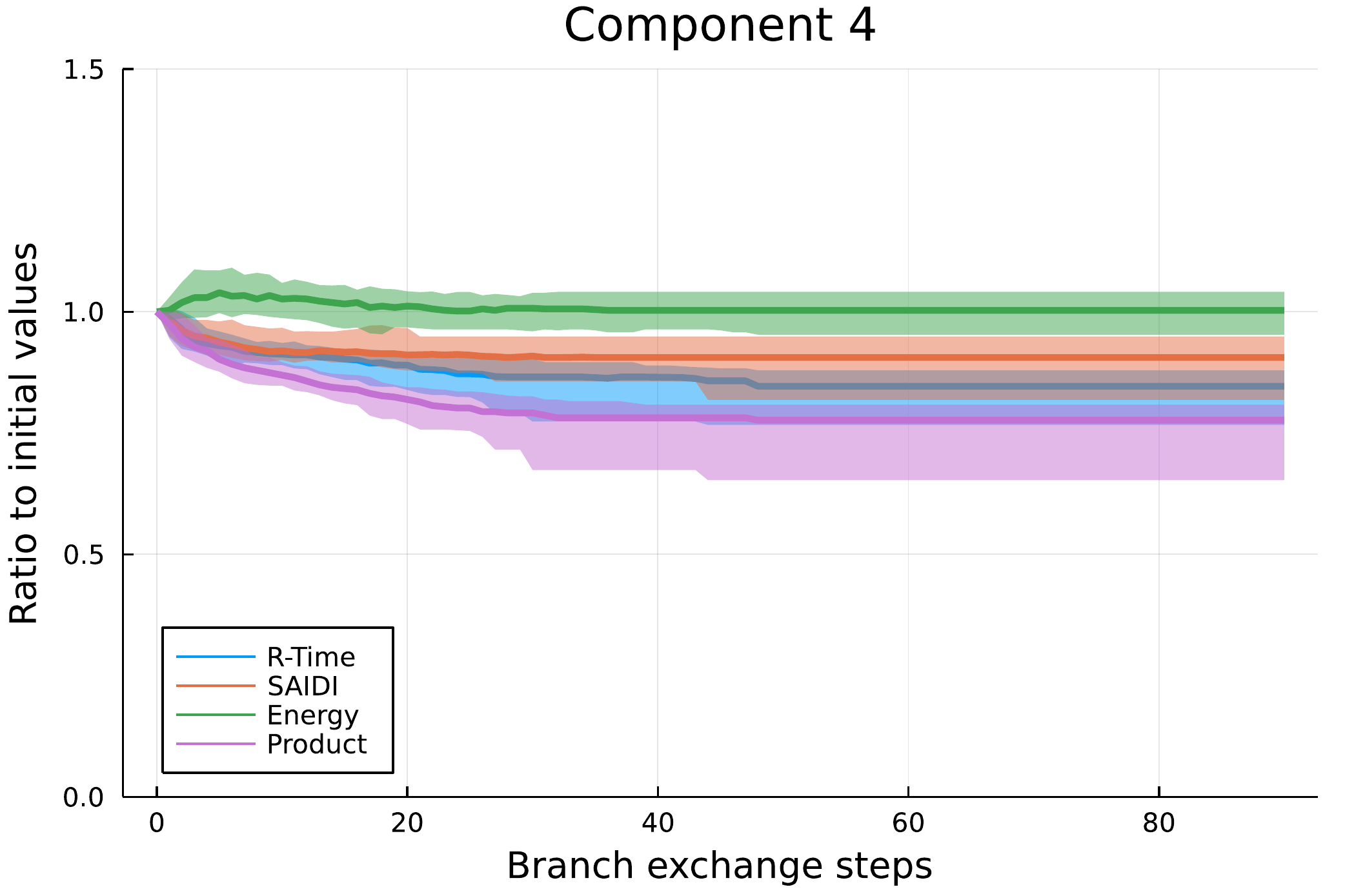}
\includegraphics[width=.49\textwidth]{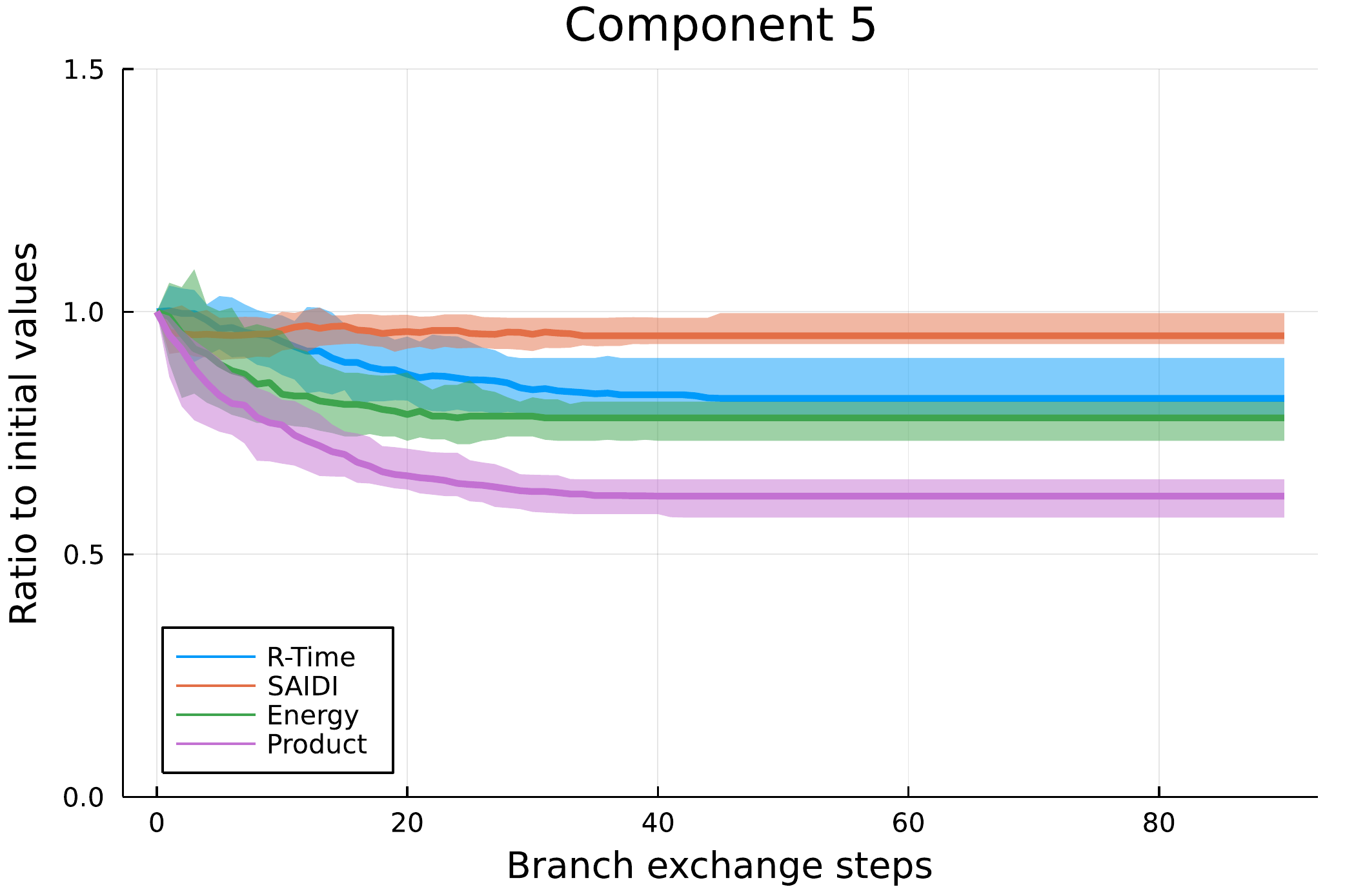}
\includegraphics[width=.49\textwidth]{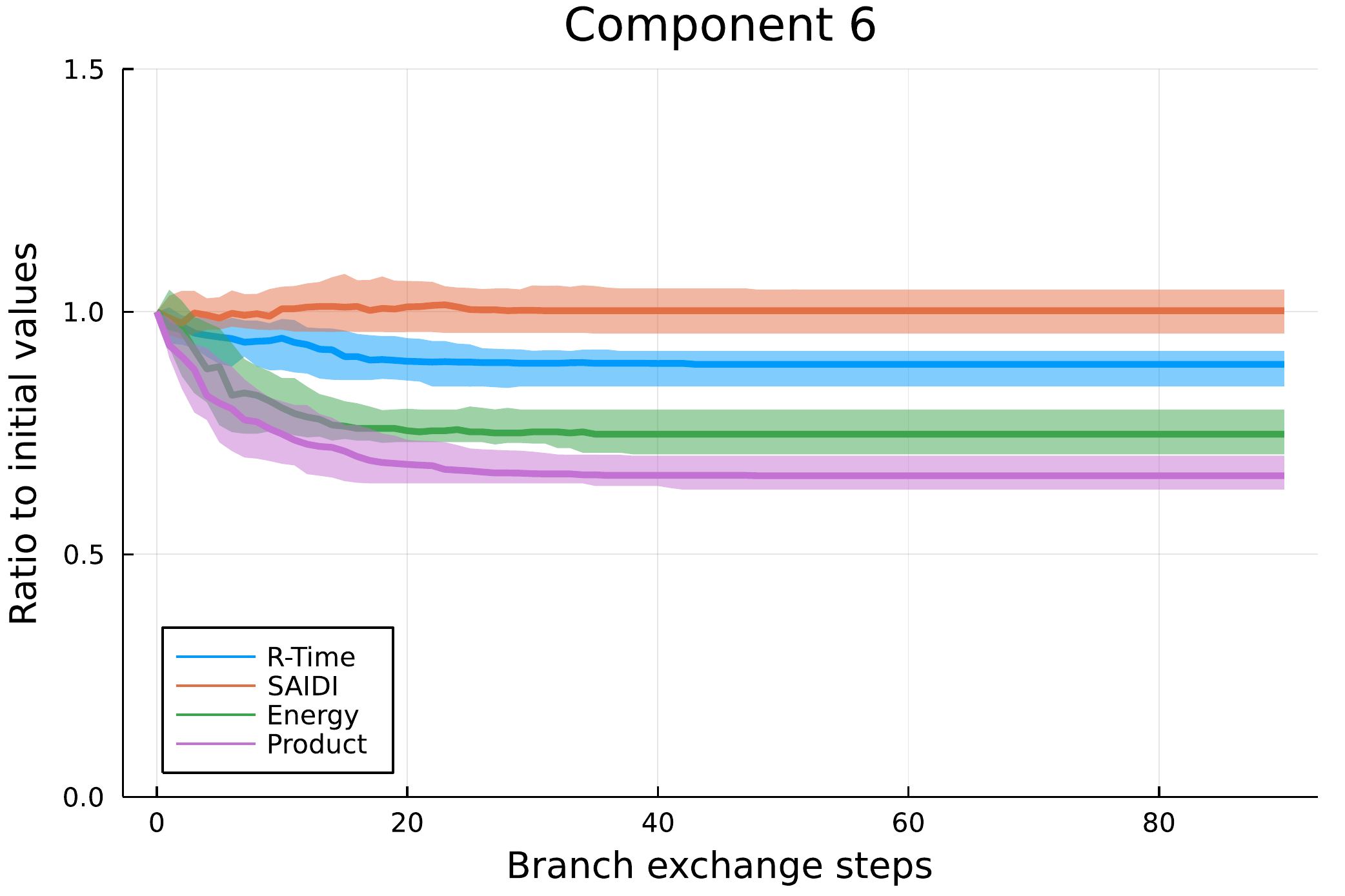}
    \caption{{\sc R-Time}, SAIDI, energy, and product objectives for Greensboro components 1-6 during the branch exchange local search. Objectives are plotted as a ratio of their value after each step to the initial values, and each plot shows the median values and percentiles (first to last deciles) for each objective, taken over 50 applications of the branch exchange local search on the component.}
\end{figure}

\clearpage
\begin{figure}[H]
    \centering
\includegraphics[width=.49\textwidth]{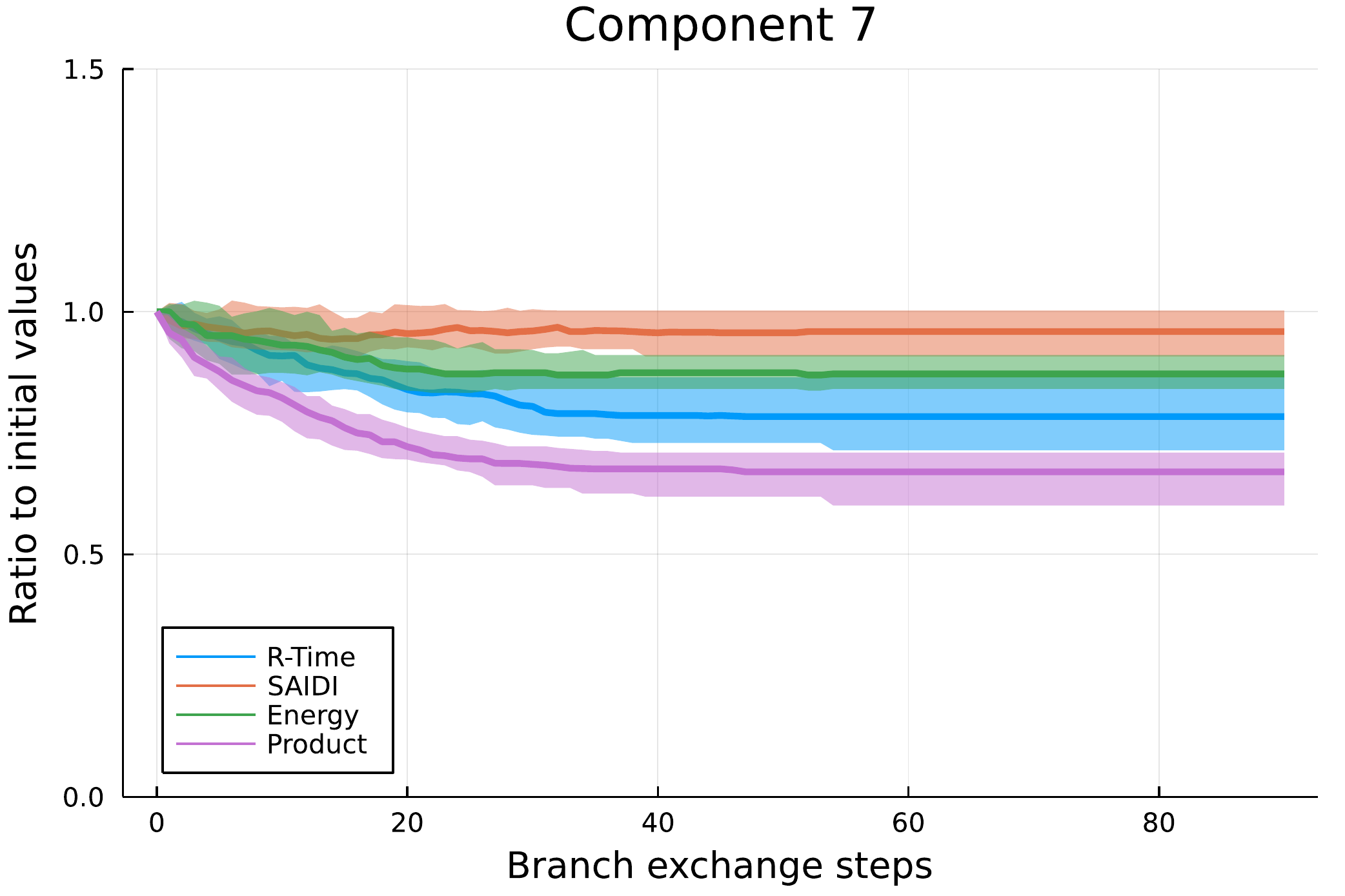}
\includegraphics[width=.49\textwidth]{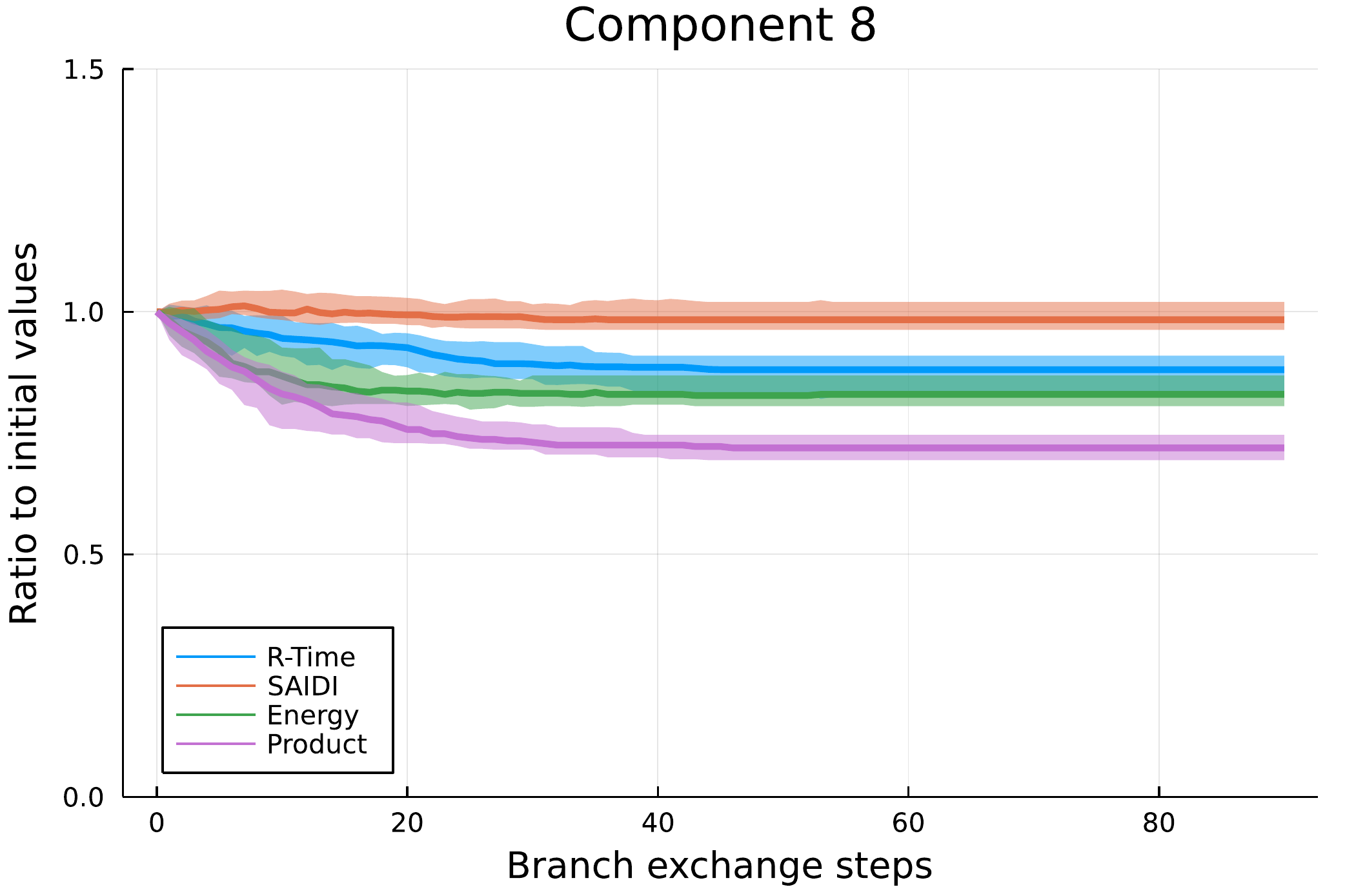}
\includegraphics[width=.49\textwidth]{Local search plots/localsearchconfidenceintervals9.pdf}
\includegraphics[width=.49\textwidth]{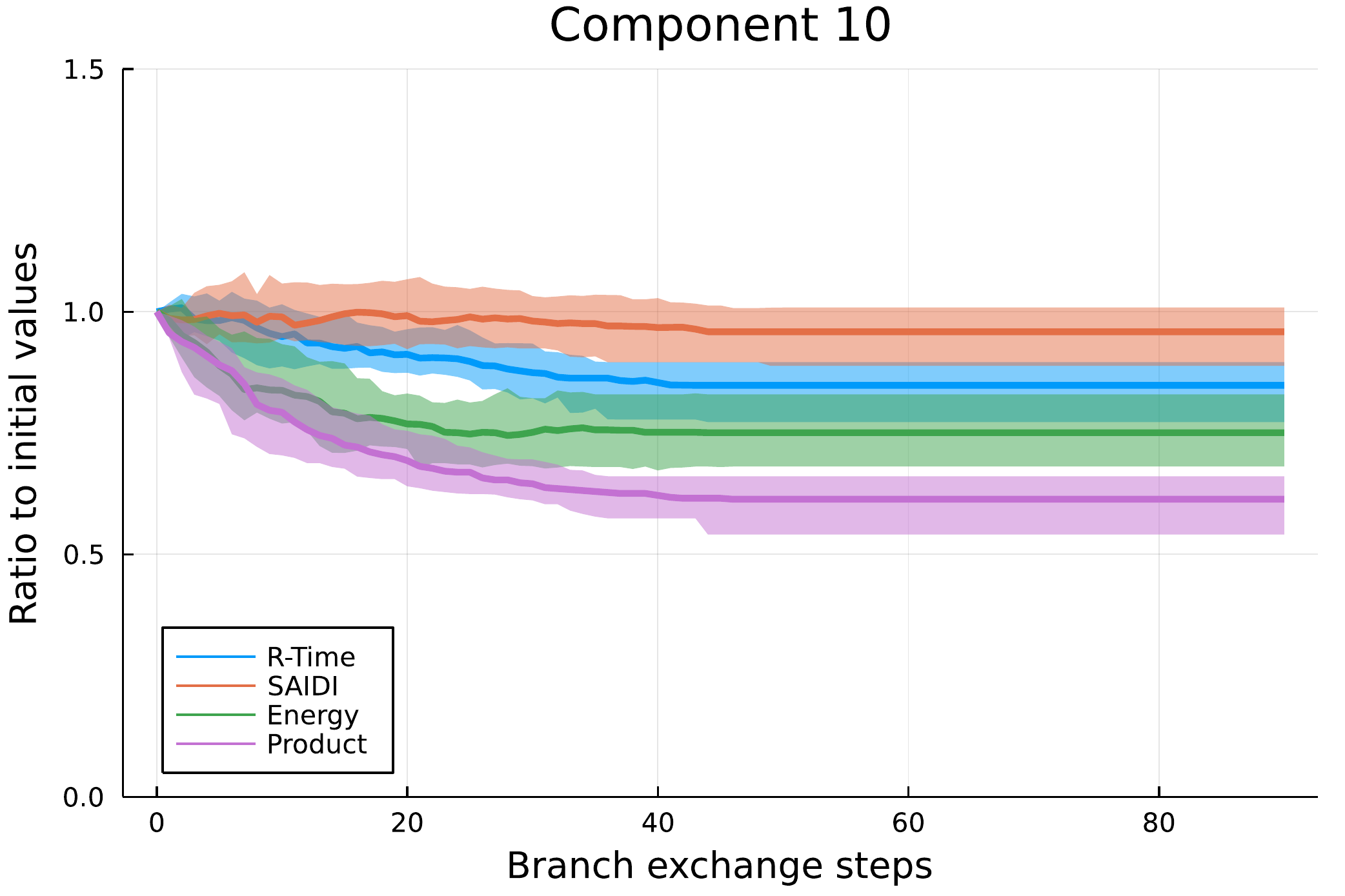}
\includegraphics[width=.49\textwidth]{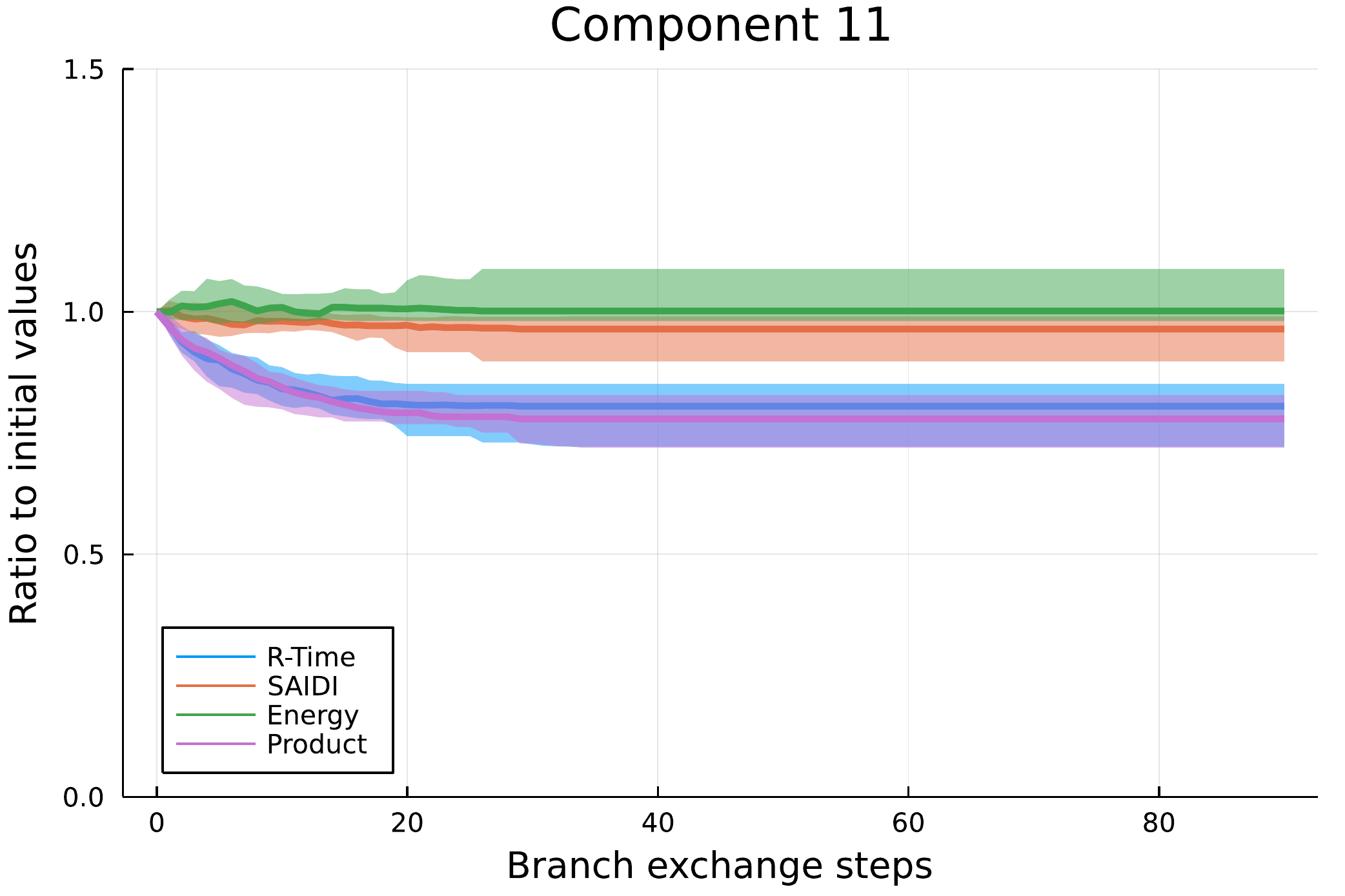}
\includegraphics[width=.49\textwidth]{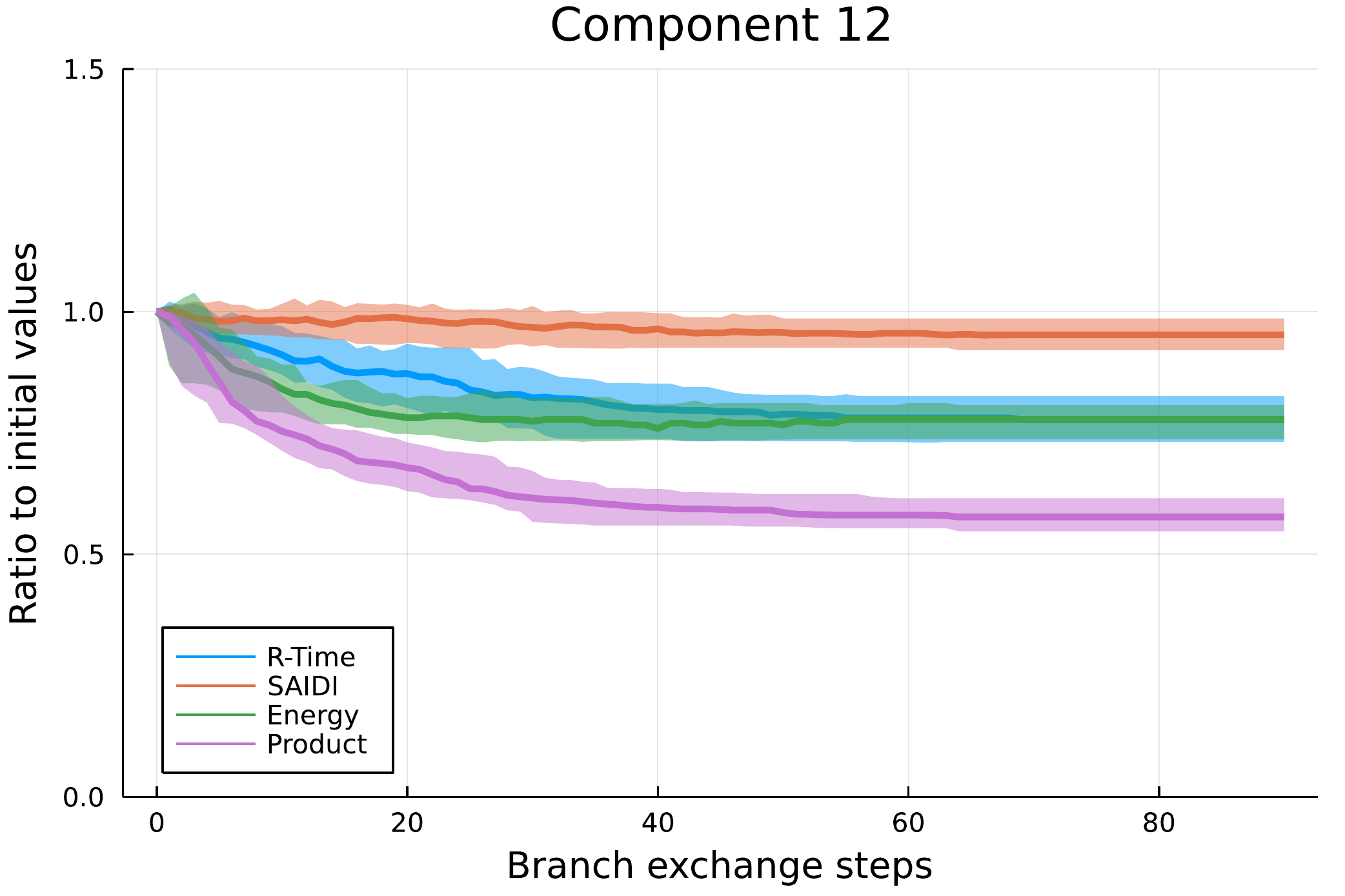}
    \caption{{\sc R-Time}, SAIDI, energy, and product objectives for Greensboro components 7-12 during the branch exchange local search. Objectives are plotted as a ratio of their value after each step to the initial values, and each plot shows the median values and percentiles (first to last deciles) for each objective, taken over 50 applications of the branch exchange local search on the component.}
\end{figure}

\clearpage

\begin{figure}[t]
    \centering
\includegraphics[width=.49\textwidth]{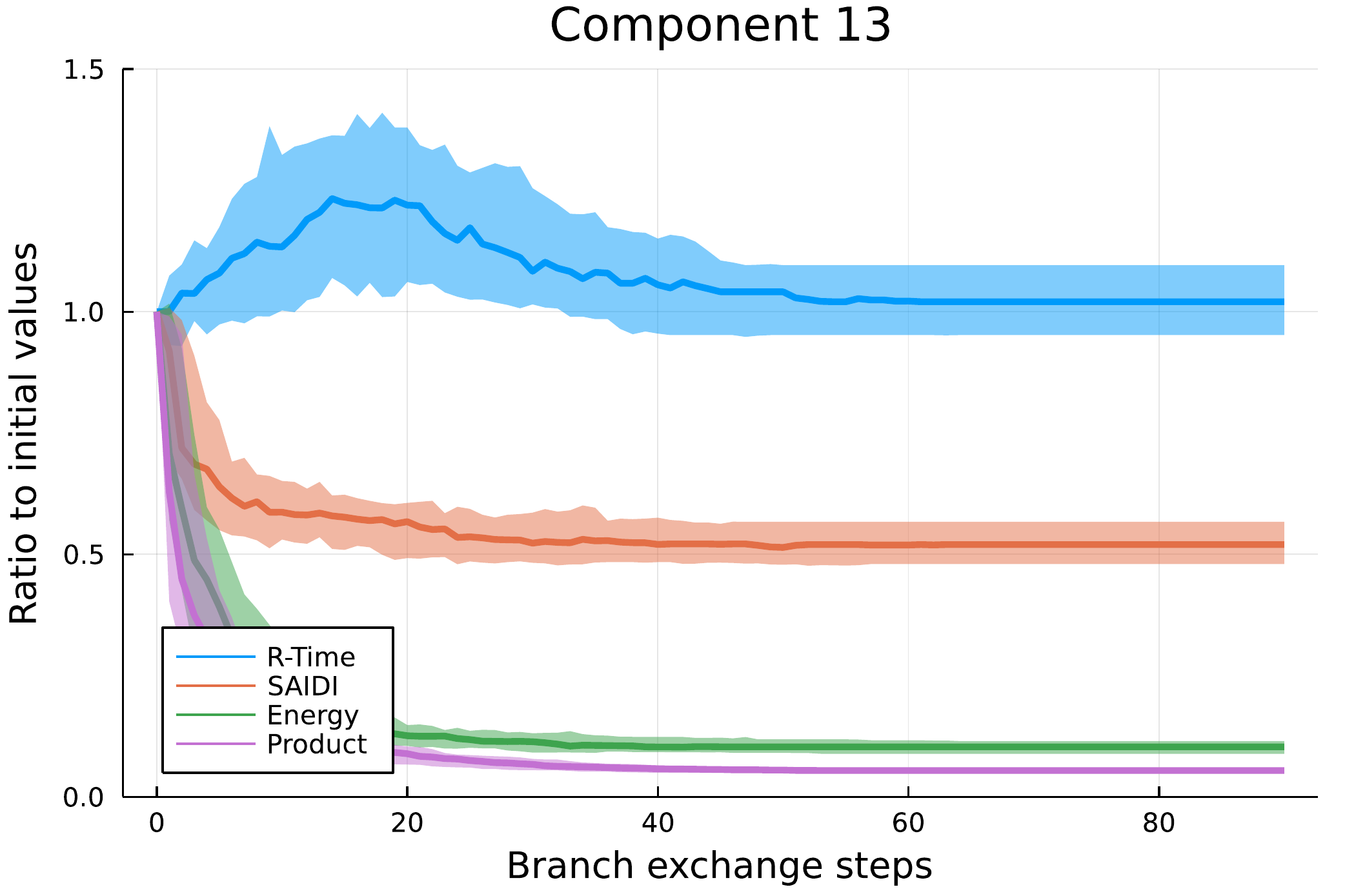}
\includegraphics[width=.49\textwidth]{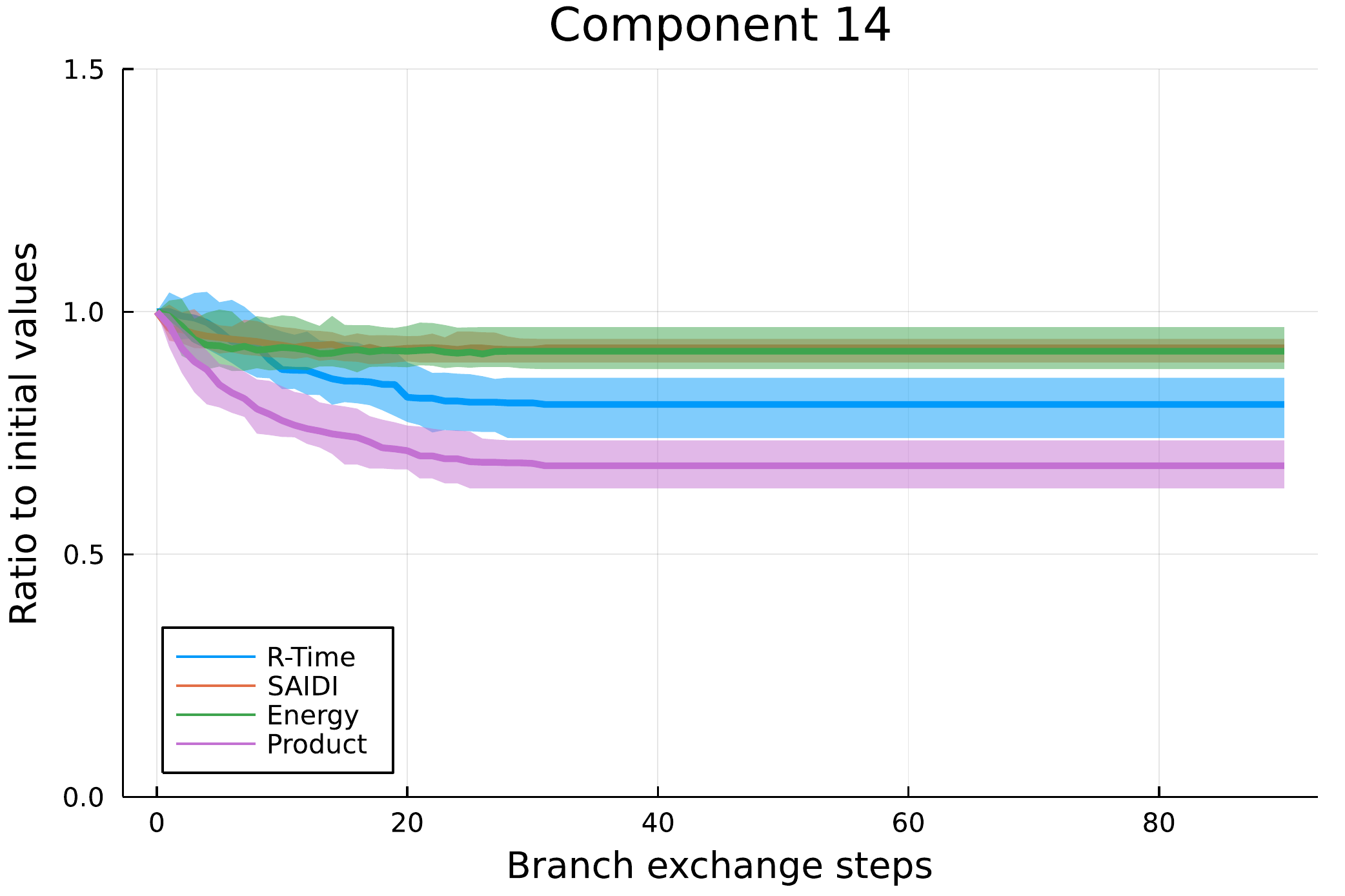}
\includegraphics[width=.49\textwidth]{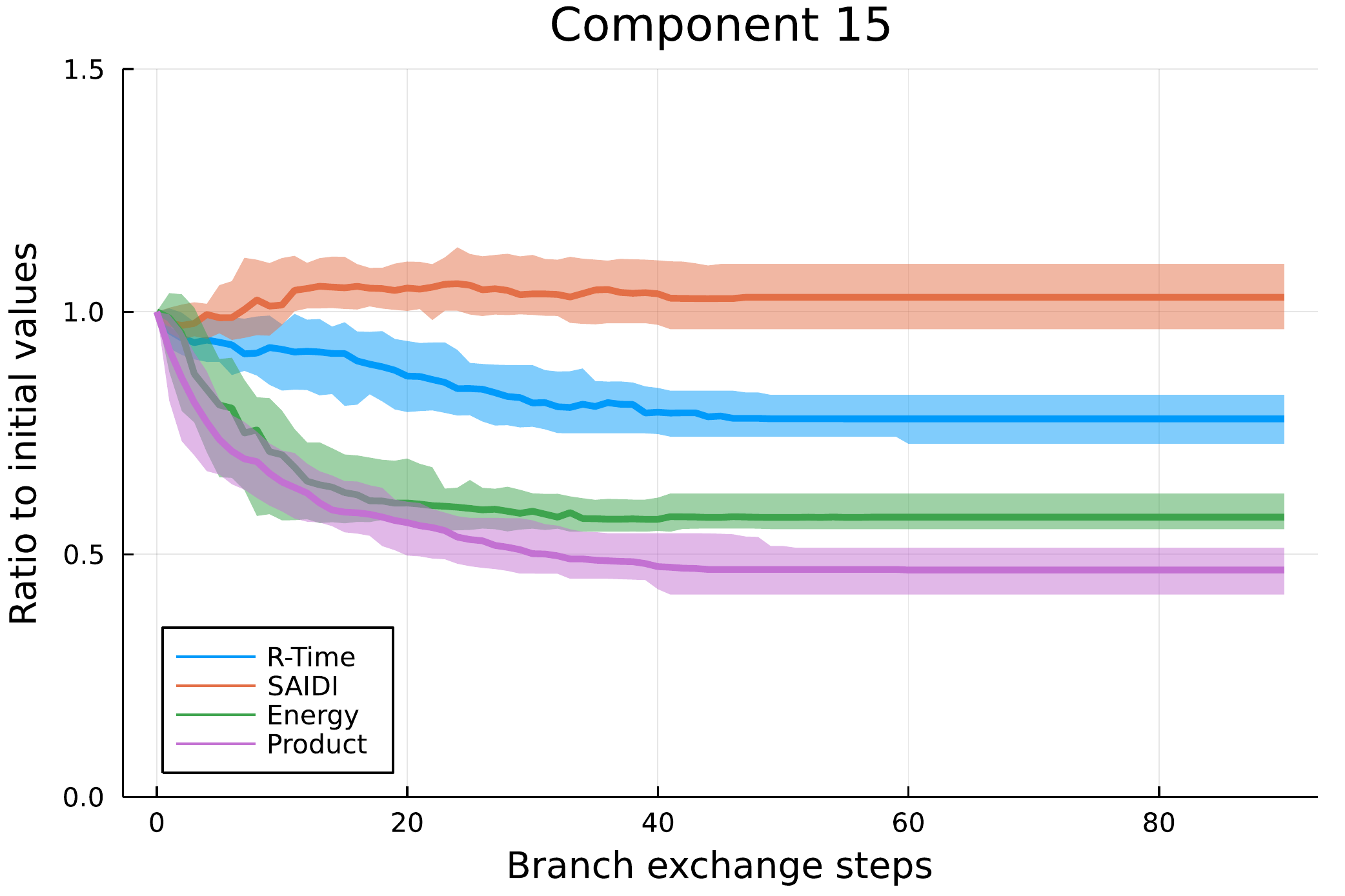}
\includegraphics[width=.49\textwidth]{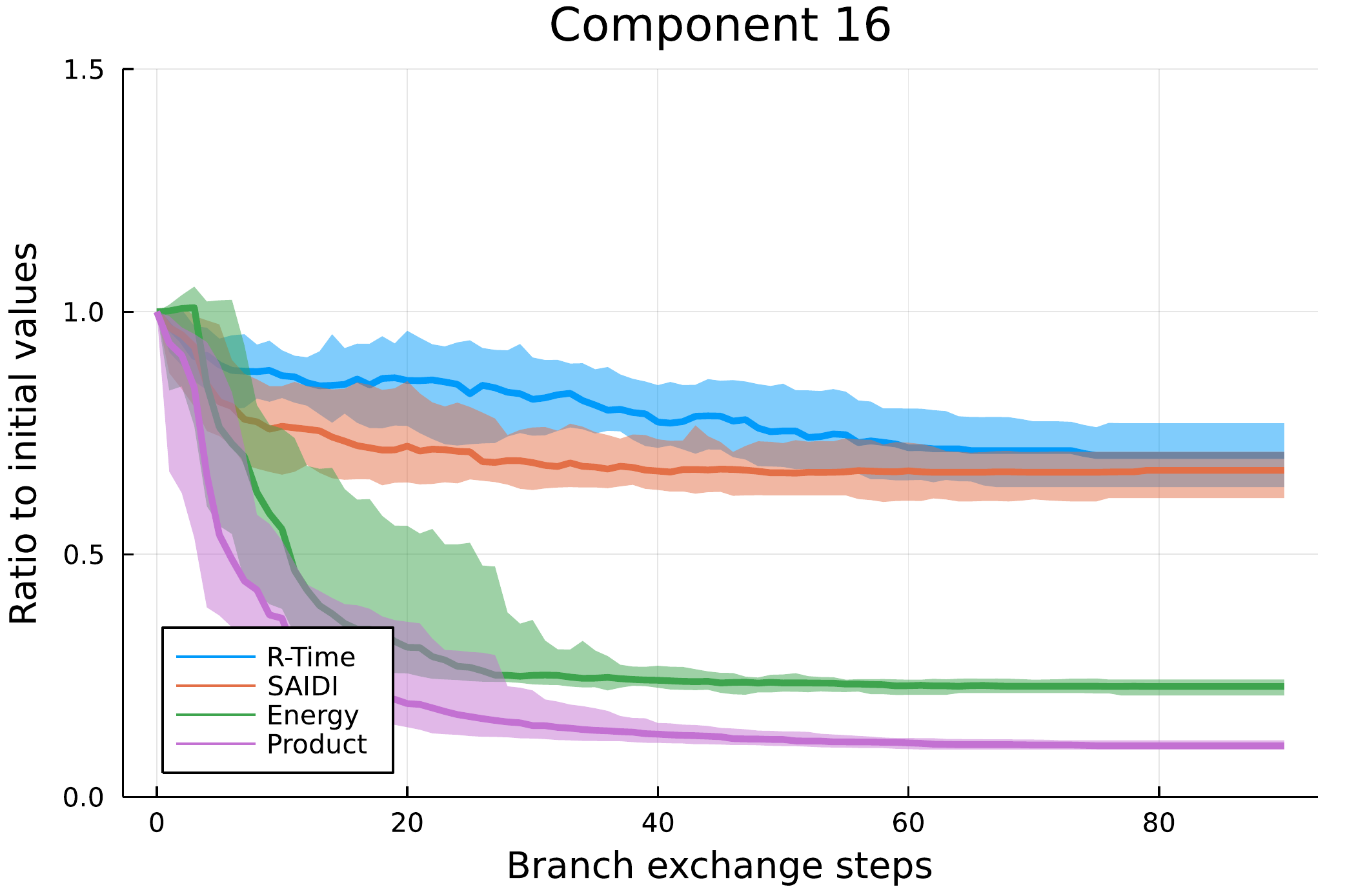}
\includegraphics[width=.49\textwidth]{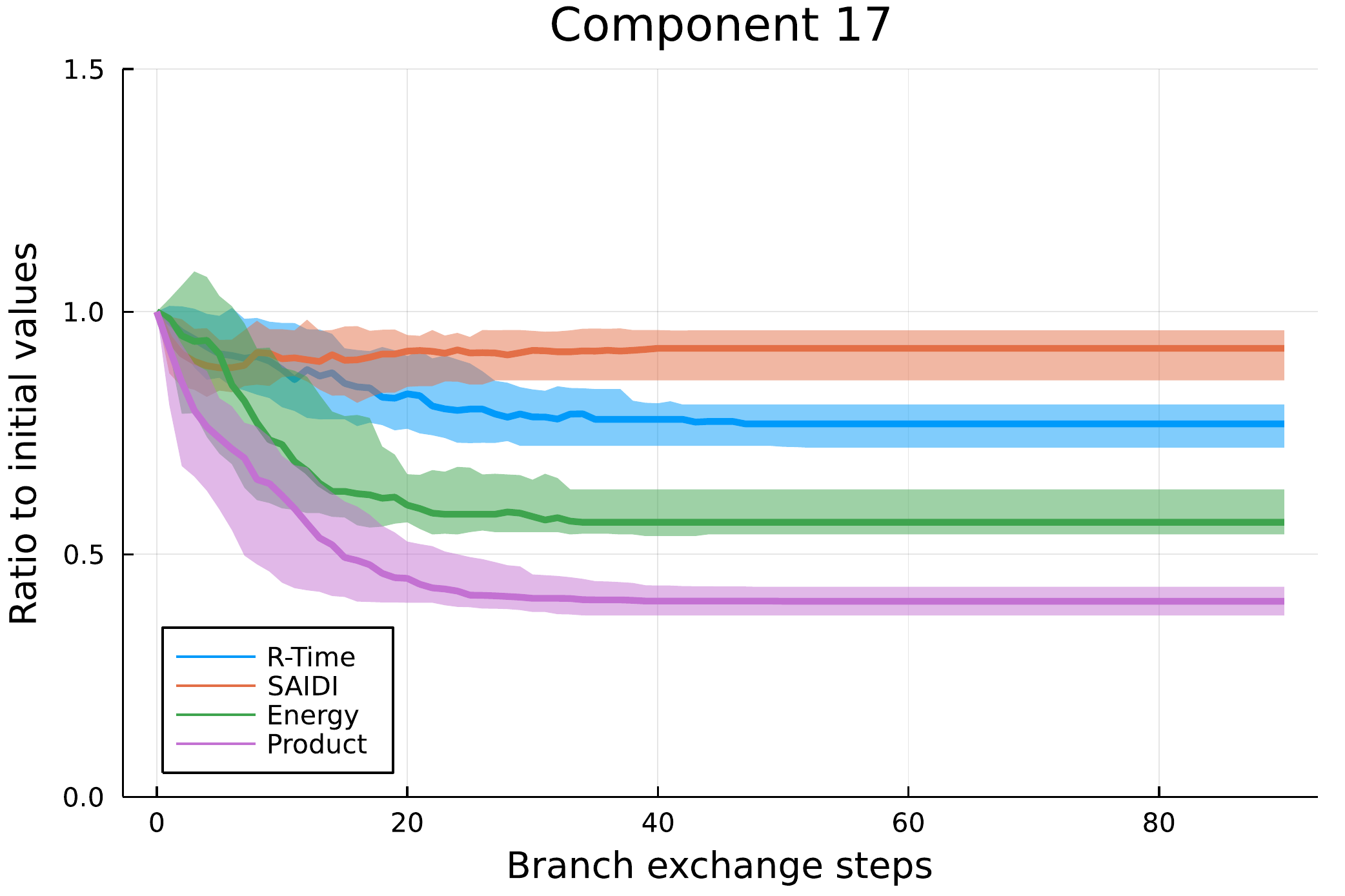}
\includegraphics[width=.49\textwidth]{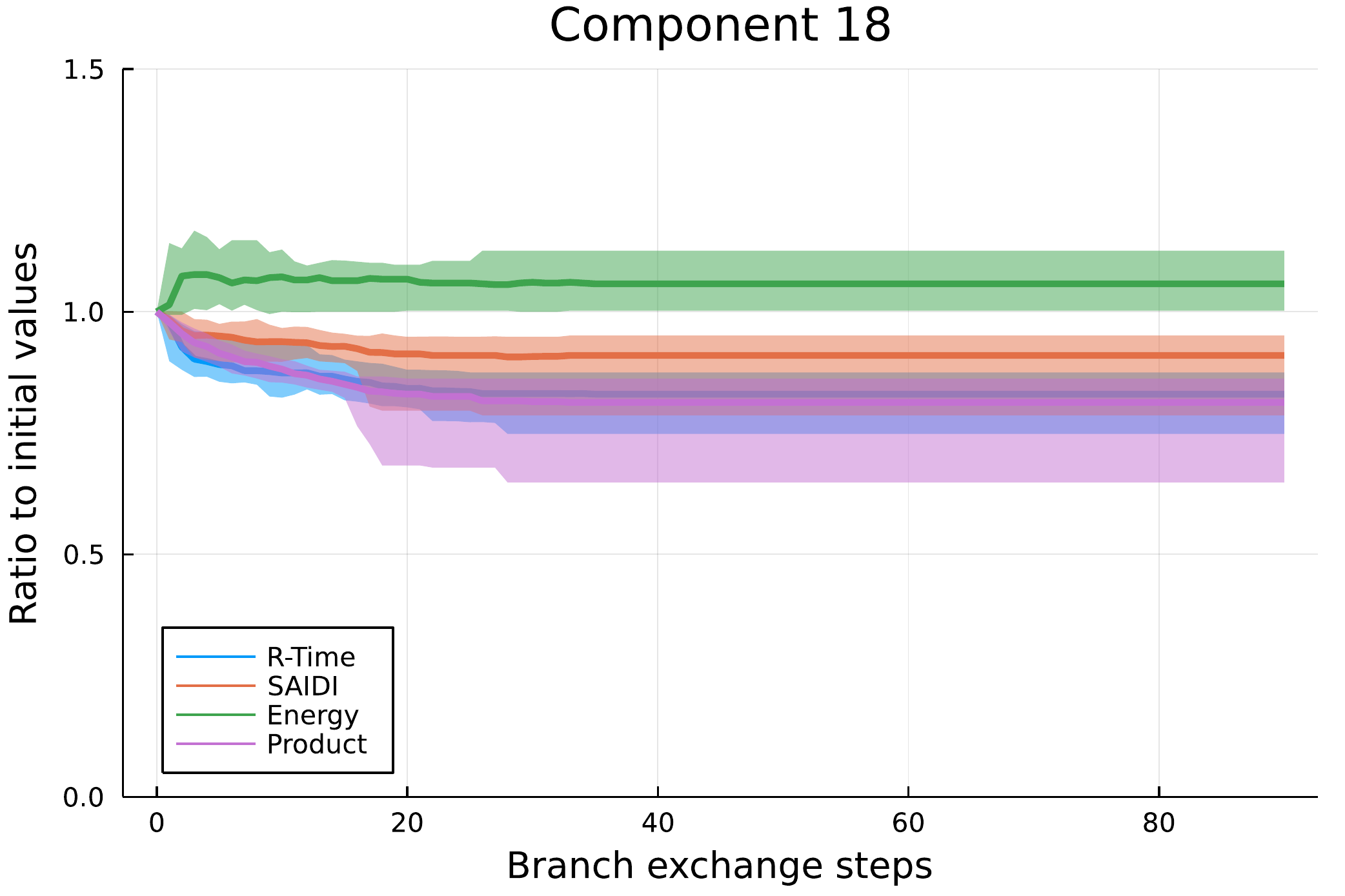}
    \caption{{\sc R-Time}, SAIDI, energy, and product objectives for Greensboro components 13-18 during the branch exchange local search. Objectives are plotted as a ratio of their value after each step to the initial values, and each plot shows the median values and percentiles (first to last deciles) for each objective, taken over 50 applications of the branch exchange local search on the component.}
    \label{fig:localsearchall}
\end{figure}

\section{Integer programming formulations}\label{sec:app_IP}
Let $G=(V,E)$ be a network with no specified initial spanning tree, and let $b(e)$ be our desired edge weights in the MRT objective $\min_{\sigma: S \longrightarrow [|S|]} \sum_{e \in T} b(e) \min_{s \text{ covers } e} \sigma(s)$ once a spanning tree $T$ is selected. By solving the following integer program, one can find both an initial spanning tree $T$ and a switch reconnection order $\sigma$ on $S=E\setminus T$ that minimizes the MRT objective:
\begin{align}
    \mathrm{minimize}\qquad & \sum_{e,t}t\cdot b(e)x_{et}\\
     \mathrm{subject\ to}\qquad & \sum_{e \subset  C} s_e \leq |C|-1, \quad \forall \text{ cuts }C\subseteq V, \label{eq:ip1_st1}\\
     \qquad & \sum_{e} s_e = n-1,\label{eq:ip1_st2}\\
     \qquad & \sum_t x_{et} = s_e, \quad \forall e,\label{eq:ip1_cover1}\\
     \qquad & x_{et} \leq \sum_{f} c_{eft} \quad \forall e,t,\label{eq:ip1_cover2}\\
     \qquad & c_{eft} \leq y_{ft}, \quad \forall e,f,t,\label{eq:ip1_cover3}\\
     \qquad & \sum_e y_{et} = 1, \quad \forall t,\label{eq:ip1_switch2}\\
     \qquad & y_{ft} \leq 1-s_f, \quad \forall f,t,\label{eq:ip1_switch1}\\
     \qquad & c_{eft} \leq \sum_{g \text{ across } C} s_g - s_e, \quad \forall e,f,t,\text{cuts $C$ crossed by $e$ but not $f$},\label{eq:ip1_cover4}\\
    \qquad & \text{all variables binary}.
\end{align}

Here constraints~\eqref{eq:ip1_st1},~\eqref{eq:ip1_st2} enforce that $T =\{e: s_e = 1\}$ is a spanning tree. Constraint~\eqref{eq:ip1_cover1} enforces that each tree edge $e$ is covered at exactly one time $t$, for which $x_{et}=1$. Constraint~\eqref{eq:ip1_cover2} enforces that $e$ is covered by a switch $f$ that has waiting time $t$ (constraints~\eqref{eq:ip1_cover3},~\eqref{eq:ip1_switch2}), is not a tree edge (constraint~\eqref{eq:ip1_switch1}), and covers $e$ (constraint~\eqref{eq:ip1_cover4}). Note that $b(e)$ should contain only $p(e)$ or other edge-specific data, as including $f(e)$ or other data that depends on $T$ will make the objective non-linear. Thus this IP can minimize MRT with {\sc R-Time} as the objective, but not with SAIDI as the objective.

This exponential-size formulation is intractable for large graphs, so we adapt the $O(mn)$-size formulation due to~\cite{martin1986sharp}. This yields an $O(m^3n)$-size formulation. Note that if the graph is planar or otherwise low genus, even smaller formulations for the spanning tree polytope and hence for this IP exist~\cite{validi2019note}.

\begin{align}
    \mathrm{minimize}\qquad & \sum_{e,t}t\cdot b(e)x_{et}\\
     \mathrm{subject\ to} \qquad & \sum_{e} s_e = n-1,\label{eq:ip2_st1}\\
     \qquad & s_e = a_{k,i,j} + a_{k,j,i}, \quad \forall e=\{i,j\}\in E, 1\leq k \leq n,\label{eq:ip2_st2}\\
     \qquad & \sum_{j} a_{k,i,j} \leq 1, \quad \forall 1\leq i,k \leq n,\label{eq:ip2_st3}\\
     \qquad & a_{k,k,j} = 0, \quad \forall 1\leq j,k\leq n\label{eq:ip2_st4}\\
     \qquad & \sum_t x_{et} = s_e, \quad \forall e,\label{eq:ip2_cover1}\\
     \qquad & x_{et} \leq \sum_{f} c_{eft} \quad \forall e,t,\label{eq:ip2_cover2}\\
     \qquad & c_{eft} \leq y_{ft}, \quad \forall e,f,t,\label{eq:ip2_cover3}\\
     \qquad & \sum_e y_{et} = 1, \quad \forall t,\label{eq:ip2_cover4}\\
    \qquad & y_{ft} \leq 1-s_f, \quad \forall f,t,\label{eq:ip2_cover5}\\
     \qquad & c_{eft} \leq r_{ef}, \quad \forall e,f,t,\label{eq:ip2_cover6}\\
     \qquad & r_{ef} \leq 1-s_e, \quad \forall e,f\label{eq:ip2_cover7}\\
     \qquad & \sum_{g} s_{e,f,g} = n-1,\quad \forall e,f\label{eq:ip2_st21}\\
     \qquad & s_{e,f,e} = 0,\quad \forall e,f\label{eq:ip2_st22}\\
     \qquad & s_{e,f,f} = 1,\quad \forall e,f\label{eq:ip2_st23}\\
     \qquad & s_{e,f,g} = s_g,\quad \forall e,f,g\neq e,f\label{eq:ip2_st24}\\
     \qquad & s_{e,f,g} = a_{e,f,k,i,j} + a_{e,f,k,j,i}, \quad \forall e,f \in E, g=\{i,j\}\in E\setminus\{ e\}, 1\leq k \leq n,\label{eq:ip2_st25}\\
     \qquad & \sum_{j} a_{e,f,k,i,j} \leq 1 + n(1-r_{ef}), \quad \forall e,f,\in E, 1\leq i,k \leq n,\label{eq:ip2_st26}\\
     \qquad & a_{e,f,k,k,j} = 0, \quad \forall e,f\in E, 1\leq j,k\leq n,\label{eq:ip2_st27}\\
     \qquad & \text{all variables binary}
\end{align}
Here constraints~\eqref{eq:ip2_st1},~\eqref{eq:ip2_st2},~\eqref{eq:ip2_st3},~\eqref{eq:ip2_st4}  enforce that $T =\{e: s_e = 1\}$ is a spanning tree, using auxiliary decision variables $a_{k,i,j}$ which are 1 if edge $e=\{i,j\}$ is directed from $i$ to $j$ in the orientation of $T$ where $k$ is the unique sink vertex. Constraints~\eqref{eq:ip2_cover1},~\eqref{eq:ip2_cover2},~\eqref{eq:ip2_cover3},~\eqref{eq:ip2_cover4},~\eqref{eq:ip2_cover5} are identical to the first formulation and encode edge waiting times based on the ordering of switches.
Variable $r_{ef}$ is 1 if switch $f$ covers edge $e$ (constraints \eqref{eq:ip2_cover6},~\eqref{eq:ip2_cover7}). To set $r_{ef}$, constraints~\eqref{eq:ip2_st21},~\eqref{eq:ip2_st22},~\eqref{eq:ip2_st23},~\eqref{eq:ip2_st24},~\eqref{eq:ip2_st25},~\eqref{eq:ip2_st26},~\eqref{eq:ip2_st27} check whether $T-e+f$ is a spanning tree. As in the constraints for the base tree $T$, $s_{e,f,g}$ is a decision variable for whether edge $g$ is in $T-e+f$ and  $a_{e,f,k,i,j}$ tracks the orientation of the edge between $i$ and $j$ when $T-e+f$ is oriented with vertex $k$ as the unique sink.

However, this formulation still lacks tractability on larger instances; for example, on contracted components of the Greensboro network it was not close to convergence after 24 hours of running. Its objective is also still limited to {\sc R-Time} and cannot incorporate varying demands or energy loss considerations.

\end{appendix}

\end{document}